\documentclass{crm-article}
\pdfoutput=1

\usepackage{citehack}
\usepackage{algorithm}
\usepackage{algorithmic}
\usepackage{hyperref}
\usepackage{multicol}
\floatname{algorithm}{Алгоритм}

\def\ag#1{{\color{black}#1}}

\newcommand{\sign}{\mathop{\mathrm{sign}}}

\begin{document}

\journalVol{10}

\journalNo{1}
\setcounter{page}{1}

\journalSection{Математические основы и численные методы моделирования}
\journalSectionEn{Mathematical modeling and numerical simulation}

\journalReceived{01.06.2016.}

\journalAccepted{01.06.2016.}


\affiliationnoref

\emailnoref

\UDC{519.85}
\title{Накопление ошибки в методе сопряжённых градиентов для вырожденных задач}
\titleeng{The error accumulation in the 
conjugate gradient method for 
degenerate problem}

\author{\firstname{А.\,Б.}~\surname{Рябцев}}
\authorfull{Антон Борисович Рябцев}
\authoreng{\firstname{A.\,B.}~\surname{Ryabtsev}}
\authorfulleng{Anton B. Ryabtsev}
\email{ryabtsev.ab@phystech.edu}
\affiliation{Национальный исследовательский университет "Московский физико-технический институт"\protect,\\ Россия, 141701, Московская область, г.Долгопрудный, Институтский пер., 9}
\affiliationeng{National Research University "Moscow Institute of Physics and Technology"\protect,\\ 9 Institutskiy per., Dolgoprudny, Moscow Region, 141701, Russia}



\begin{abstract}
В данной работе рассматривается метод сопряжённых градиентов при решении задачи минимизации квадратичной функции с аддитивным шумом в градиенте. Были рассмотрены три концепции шума: враждебный шум в линейном члене, стохастический шум в линейном члене и шум в квадратичном члене, а также комбинации первого и второго с последним. Экспериментально получено, что накопление ошибки отсутствует для любой из рассмотренных концепций, что отличается от фольклорного мнения, что, как и в ускоренных методах, накопление ошибки должно иметь место. В работе приведена мотивировка того, почему ошибка может и не накапливаться. Также экспериментально исследовалась зависимость ошибки решения как от величины (масштаба) шума, так и от размера решения при использовании метода сопряжённых градиентов. Предложены и проверены гипотезы о зависимости ошибки в решении от масштаба шума и размера (2-нормы) решения для всех рассмотренных концепций. Оказалось, что ошибка в решении (по функции) линейно зависит от масштаба шума. В работе приведены графики, иллюстрирующие каждое отдельное исследование, а также детальное описание численных экспериментов, включающее в себя изложение способов зашумления как вектора, так и матрицы.
\end{abstract}


\begin{abstracteng}
In this paper, we consider the conjugate gradient method for solving the problem of minimizing a quadratic function with additive noise in the gradient. Three concepts of noise were considered: antagonistic noise in the linear term, stochastic noise in the linear term and noise in the quadratic term, as well as combinations of the first and second with the last. It was experimentally obtained that error accumulation is absent for any of the considered concepts, which differs from the folklore opinion that, as in accelerated methods, error accumulation must take place. The paper gives motivation for why the error may not accumulate. The dependence of the solution error both on the magnitude (scale) of the noise and on the size of the solution using the conjugate gradient method was also experimentally investigated. Hypotheses about the dependence of the error in the solution on the noise scale and the size (2-norm) of the solution are proposed and tested for all the concepts considered. It turned out that the error in the solution (by function) linearly depends on the noise scale. The work contains graphs illustrating each individual study, as well as a detailed description of numerical experiments, which includes an account of the methods of noise of both the vector and the matrix.
\end{abstracteng}
\keywordeng{conjugate gradient method}
\keywordeng{degenerate problem}
\keywordeng{noisy oracle.}

\maketitle

\paragraph{Введение}
Во многих приложениях часто  необходимо решать систему линейных алгебраических уравнений:
$$Ax=b.$$
Если точное решение не требуется, а нужно лишь найти какое-то приближение, то данную задачу можно свести к задаче минимизации квадратичной функции (считаем матрицу $A$ симметричной и неотрицательно определенной):
$$ \;\; f(x)=\frac{1}{2}\langle Ax,x\rangle-\langle b,x\rangle\rightarrow \min_{x\in\mathbb{R}^n}$$
и решать её известными методами численной оптимизации, такими как, например, ускоренный метод Нестерова \cite{nesterov10}, или методом сопряжённых градиентов \cite{Gasnikov17}. 

Многие задачи, приходящие из реальных приложений, оказываются вырожденными (наименьшие собственные значения матрицы $A$ равны нулю или близки к нулю), см., например, \cite{kabanikhin2012}. В данной статье также будем рассматривать вырожденный случай.

Известно\footnote{Подробнее в разделе \nameref{related works}.}, что ускоренные методы могут оказаться неустойчивыми к неточностям в градиентах, что приводит к накоплению ошибки с ростом числа итераций. Насколько нам известно, метод сопряжённых градиентов до сих пор не был в теоретическом плане исследован в данном ключе, хотя, безусловно, также является ускоренным. По-видимому, причина связана с наличием отрицательных результатов о сходимости методов с одномерным поиском. Так, например, как следует из \cite{poljak81} метод наискорейшего спуска со сколь угодно малой аддитивной неточностью в градиенте может в итоге расходиться. Однако, в данной работе показано, что типично такие ситуации для метода сопряжённых градиентов на квадратичных задачах не возникают и шум не накапливается.


Стоит отметить, что задача минимизации положительно определенной квадратичной формы является классической задачей выпуклой оптимизации. Исследование данного семейства задач может дать представление о сходимости (хотя бы в окрестности решения) различных методов в задачах выпуклой оптимизации.  Метод сопряжённых градиентов гарантированно находит точное решение этой задачи за $N = n$ итераций, где $n$~---~размер задачи. Это свойство является отличительной особенностью методов сопряжённых градиентов от всевозможных
обобщений. Однако, трудности появляются при решении вырожденных (некорректных) задач квадратичной оптимизации.
В этом случае, число требуемых итераций может быть близко к $n$. С учетом того, что стоимость итерации $O(n^2)$, получается общая трудоемкость $ O(n^3)$, что можно получить и более простыми алгоритмами, например, методом Гаусса.

В данной статье мы рассмотрим простейшую задачу:
$$Ax = b,$$
в условиях отсутствия точных значений матрицы $A$ и(или) вектора $b$. Нам доступны матрица $\tilde{A}$ вместо исходной $A$ или вектор $\tilde{b}$ вместо исходного $b$:
$$
\|\tilde{A}-A\|_2 \le \delta_A, \;\;\;\|\tilde{b}-b\|_2 \le \delta_b.
$$

\smallskip

Оба случая рассмотрены как отдельно, так и вместе. Поставленная задача сводится к задаче минимизации квадратичной формы:

$$1) \;\; f(x)=\frac{1}{2}\langle Ax,x\rangle-\langle\tilde{b},x\rangle\rightarrow \min_{x\in\mathbb{R}^n},$$

$$2) \;\; f(x)=\frac{1}{2}\langle \tilde{A}x,x\rangle-\langle b,x\rangle\rightarrow \min_{x\in\mathbb{R}^n},$$

$$3) \;\; f(x)=\frac{1}{2}\langle \tilde{A}x,x\rangle-\langle \tilde{b},x\rangle\rightarrow \min_{x\in\mathbb{R}^n}.$$

\medskip

\paragraph{Связанные исследования}
\label{related works}

Известно \cite{devolder2013}, \cite{dvinskikh2019}, что при использовании неточного градиента $\tilde{\nabla}f(x)$, удовлетворяющего для всех $x$, $y$
$$
f(x) + \langle\tilde{\nabla}f(x), y-x\rangle - \delta_1 \le f(y) \le
f(x) + \langle\tilde{\nabla}f(x), y-x\rangle + \dfrac{L}{2}\|y-x\|_2^2 + \delta_2
$$
справедлива оценка 
$$
f(x^N) - f(x^*) = O\Big(\dfrac{LR^2}{N^p} + \delta_1 + N^{p-1}\delta_2\Big).
$$
Здесь $p\in\{1, 2\}$ для неускоренных и ускоренных методов соответсвенно, а $R=\|x_0-x^*\|_2$.

Для случая аддитивного шума в градиенте 
$$
\|\tilde{\nabla}f(x) - \nabla f(x)\|_2
\ag{\le} 
\delta
$$
в цикле работ А.С. Немировского \cite{ nemirovskiy84}, \cite{nemirovskiy86}, \cite{nemirovskiy92}  были получены интересные результаты о регуляризующих свойствах метода сопряжённых градиентов для вырожденных (некорректных) задач квадратичной оптимизации. Вырожденной будем называть задачу выпуклой оптимизации, для которой отношение максимального и минимального собственного значения функционала (обусловленность задачи) много больше квадрата
размерности пространства, в котором происходит оптимизация: $\dfrac{L}{\mu} \gg n^2$, и не меньше величины обратной к относительной точности, с
которой требуется решить задачу. Например, к такому классу задач относятся задача минимизации квадратичной формы, заданной матрицей с набором собственных чисел как на рис. 1.
\smallskip
\begin{figure}[ht]
	\centering
	\includegraphics[scale=0.3]{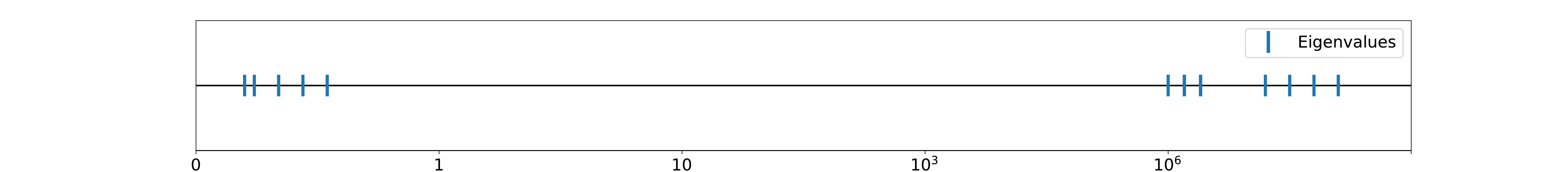}
	\caption{Спектр плохообусловленной матрицы.}
	\label{fig:fig1}
\end{figure}
\smallskip

Многие задачи, приходящие из реальных приложений, оказываются вырожденными, см., например, \cite{kabanikhin2012}.
Строить сходящиеся по аргументу алгоритмы для таких задач в общем
случае оказывается невозможным. Решение задачи оказывается неустойчивым к неточностям в данных. Для возможности корректного восстановления решения требуются дополнительные предположения (истокопредставимости). Здесь мы ограничимся простейшей задачей:
$$
Ax=b,
$$
в которой не доступны точные значения $A$ и $b$ , а доступны только $\tilde{A}$ и $\tilde{b}$, где
$$
\|\tilde{A} - A\|_2 \le \delta_A, \;\; \|\tilde{b} - b\|_2 \le \delta_b,
$$
где $\bigl\|C\bigl\|_2 = \sqrt{\lambda_{max}(C^TC)}$.
По поставленной задаче можно построить следующие задачи оптимизации
$$
1) \;\; f_1(x) = \dfrac{1}{2}\Bigl\|\tilde{A}x - \tilde{b}\Bigr\|_2^2 \rightarrow \min_{x\in \mathbb{R}^n},
$$

$$
2) \;\; f_2(x) = \dfrac{1}{2}\langle\tilde{A}x, x\rangle - \langle \tilde{b}, x \rangle \rightarrow \min_{x\in \mathbb{R}^n}\; \text{(если }A^T=A,\; \tilde{A}^T = \tilde{A}).
$$

\smallskip

Введём индекс $\tau \in \{1, 2\}$, который будет отвечать рассматриваемому случаю. В работе \cite{nemirovskiy86} было показано, что в случае, когда выполняется условие истокопредставимости
$$
x_*=(A^TA)^{\sigma/2}y_*, \;\; \bigl\|y_*\bigr\|_2 \le R_{\sigma}, \;\; Ax_* = b,
$$
метод сопряжённых градиентов с критерием останова вида
$$
\Bigl\|\tilde{A}x^N-\tilde{b}\Bigr\|_2 \le 2\Bigl(\delta_A \Bigl\|x^N\Bigr\|_2 + \delta_b\Bigr),
$$
стартующий с $x_0 = 0$, сходится для соответствующей задачи $\tau \in \{1, 2\}$ следующим образом

\smallskip

$$
\omega_N^2 = \Bigl\|\tilde{A}x^N-\tilde{b}\Bigr\|_2^2 = O\Biggl(\dfrac{\tilde{L}^{2(1+\sigma)}R_{\sigma}^2}{N^{2\tau(1+\sigma)}} + \omega_*^2\Biggl), \; \omega_* = \tilde{L}^{\sigma}R_{\sigma}\delta_A + \delta_b,
$$
\smallskip
где $\tilde{L} = \max \Bigl\{\bigl\|A\bigr\|_2, \bigl\|\tilde{A}\bigr\|_2\Bigr\}$, причем до выполнения критерия останова
\smallskip
$$
\Bigl\|\tilde{A}x^N-\tilde{b}\Bigr\|_2 \le 2\Bigl(\delta_A \Bigl\|x^N\Bigr\|_2 + \delta_b\Bigr)
$$
\smallskip
при $\theta + 2\sigma > 0$, $\theta \in [0, 2]$ справедлива следующая оценка
\smallskip
$$
\nu_{\theta, N}^2 = \Bigl\|\bigl(A^TA\bigr)^{\theta/4}\bigl(x^N-x_*\bigr)\Bigr\|_2^2 = O\Bigl(R_{\sigma}^{(2-\theta)/(1+\sigma)}w_N^{(\theta+2\sigma)/(1+\sigma)}\Bigr),
$$
\smallskip
$$
\bigl\|x^N\bigr\|_2 = O\bigl(\bigl\|x_*\bigr\|_2\bigr).
$$
\smallskip

Обратим внимание, что в $\nu_{\theta, N}^2$ стоит настоящая (незашумленная) матрица $A$. Приведенные выше результаты являются точными и не могут быть улучшены за счет использовании других методов. Причем не могут быть улучшены как в части скорости сходимости, так и в части достижимой точности $O(\omega_*)$.  Удивительно здесь, в частности, то, что метод сопряжённых градиентов, безусловно, можно относить к классу ускоренных (оптимальных) методов, для которых известно, что в общем случае, неточность в вычислении градиента линейно накапливается с ростом номера итерации \cite{devolder2013}. Однако приведенный выше результат свидетельствует об отсутствии накопления неточностей, что соответствует неускоренным методам.

\ag{По-видимому, связано это со спецификой шума -- шум в градиенте аддитивный. В работе \cite{dAspremont} рассмотрена задача $f(x)\rightarrow \min_{x\in Q}$ с компактным $Q$.
В этой работе показано, что для ускоренных методов аддитивный шум в градиенте не накапливается по мере роста итераций. Более общий результат (не требующий компактности $Q$) приведен в работе \cite{dvinskikh2019}.}

\ag{Естественно, появляется гипотеза, что и для метода сопряжённых градиентов подобно тому, что имеет место для ускоренных градиентных методов, не будет наблюдаться накопление неточностей по мере роста итераций, если шум в градиенте аддитивен. Об этом, в частности, говорят результаты, приведенные выше, в которых шум аддитивен и не менялся с номером итерации. }
\ag{Более точно} в данной работе проверялась 
гипотеза,
что для случая $\|\tilde{\nabla}f(x) - \nabla f(x)\|_2 \le \delta$ \ag{для метода сопряжённых градиентов} верна оценка
$$
f(x^N) - f(x^*) = O\Big(\dfrac{LR^2}{N^2} + \delta R\Big).
$$
Причём для шума в векторе $b$ отдельно рассматривались две концепции: враждебный шум, зависящий от направления градиента в данный момент, и случайный, прибавляющийся или вычитающийся из вектора $b$ на каждой итерации с вероятностью $\dfrac{1}{2}$.

Если предполагать, что последовательность, генерируемая методом, ограничена, то шум в матрице можно рассматривать тоже как аддитивный шум в градиенте. Для данной концепции проверена гипотеза

\smallskip

$$
\big|f(x^*_{noisy}) - f(x^*)\big| = O\big(\delta R^2\big).
$$

\bigskip

Здесь $f(x^*_{noisy})$~---~это то, к чему стремится\footnote{Можно считать, что $f(x^k)$ стремится к $f(x^*_{noisy})$ в Чезаровском смысле (то есть в смысле средних арифметических).} функционал незашумлённой задачи на последовательности, сгенерированной сопряжёнными градиентами с неточным градиентом.

\vspace{0.2cm}

\paragraph{Концепции шума в векторе $b$}
\label{conceptions}

\vspace{0.1cm}

В работе проведены эксперименты для двух концепций шума в векторе $b$: 1) враждебного и 2) случайного.

\vspace{0.2cm}

\begin{enumerate}
    \item $\delta^k = \delta sign(Ax^k-b)$
    
    \vspace{0.1cm}
    
    \item $\delta^k =\begin{cases}\hspace{0.29cm}\delta,\; \text{с вероятностью}\;\; \frac{1}{2} \\-\delta,\; \text{с вероятностью}\;\; \frac{1}{2}\end{cases}$
\end{enumerate}

\vspace{0.2cm}

Суть враждебного метода в отодвигании на каждой итерации точки оптимума в направлении антиградиента, что как бы создаёт "убегание"\protect, которое препятствует своевременному нахождению решения.

Отметим, что при исследованиях поведения метода в случае зашумления матрицы, шум на каждой итерации вносился согласно идее случайного шума для вектора, то есть сгенерированная матрица шума прибавлялась к исходной равновероятно то со знаком плюс, то со знаком минус.

Подробнее в разделе \nameref{supplement materials}.

\medskip
\medskip

\paragraph{Метод сопряжённых градиентов}

Для экспериментов в работе использовался метод сопряжённых градиентов, а именно~---~алгоритм Флетчера-Ривса (алгоритм \ref{Alg}).

\begin{algorithm}[!h]
\caption{{Метод сопряжённых градиентов (алгоритм Флетчера-Ривса).}}
\label{Alg}
\hspace*{\algorithmicindent} \textbf{Input: } $x_0$~---~начальная точка, $d_0=-\nabla f(x_0)$.
\begin{algorithmic}[1]
\STATE $i := 0$
\WHILE{$\|d_i\|_2 \ge \varepsilon$}
\STATE \label{loop_state} Вычислить $\alpha_i$, минимизирующее $f(x_i+\alpha_i d_i)$ по формуле:
\begin{equation}\label{exitLDL_G}
\alpha_i = -\dfrac{d_i^T(Ax_i+b)}{d_i^TAd_i}
\end{equation}
Сделать шаг:
\begin{equation}\label{exitLDL_G}
x_{i+1} = x_i + \alpha_i d_i
\end{equation}
Обновить направление:
\begin{equation}\label{exitLDL_G}
d_{i+1} = -\nabla f(x_{i+1}) + \beta_i d_i,
\end{equation}
где $\beta_i$ вычисляется по формуле:
\begin{equation}\label{exitLDL_G}
\beta_i = \dfrac{\nabla f(x_{i+1})^T A d_i}{d_i^T A d_i}
\end{equation}
\ENDWHILE
\end{algorithmic}
\hspace*{\algorithmicindent} \textbf{Output: } $x_N$
\end{algorithm}

\vspace{-0.5cm}

\paragraph{Численные эксперименты}

\vspace{-0.4cm}

Для графиков, иллюстрирующих исследования выхода траекторий функции на асимптоту, были использованы параметры из таблицы \ref{table:1}.

\begin{table}[h!]
\centering
\begin{tabular}{|l|c|c|c|c|}
\hline
                        & n      & $\delta_A$            & $\delta_b$ & $R$            \\ \hline
Враждебный шум          & $10^4$ &                       & $0.1$      & $\approx 7000$ \\ \hline
Стохастический шум      & $10^4$ &                       & $0.1$      & $\approx 7000$ \\ \hline
Шум в матрице           & $10^3$ & $\{0.0025;\; 0.005\}$ &            & $\approx 2000$ \\ \hline
Шум в матрице и столбце & $10^3$ & $0.005$               & $0.1$      & $\approx 2000$ \\ \hline
\end{tabular}
\caption{Параметры задачи при исследовании выхода траекторий функции на асимптоту. Обозначения: размерность задачи $n$, размер шума в матрице и в векторе соответственно $\delta_A$ и $\delta_b$, размер решения $R = \|x^* - x_0\|$.}
\label{table:1}
\end{table}


Графики, относящиеся к исследованию зависимости ошибки от $\delta$ при решении зашумлённой задачи, соответствуют экспериментам с параметрами из таблицы \ref{table:2}.

\begin{table}[h!]
\centering
\begin{tabular}{|l|c|c|c|c|}
\hline
                   & n      & $\delta_A$    & $\delta_b$   & $R$            \\ \hline
Враждебный шум     & $10^4$ &               & $[0;\; 0.1]$ & $\approx 2000$ \\ \hline
Стохастический шум & $10^4$ &               & $[0;\; 0.1]$ & $\approx 2000$ \\ \hline
Шум в матрице      & $10^3$ & $[0;\; 0.01]$ &              & $\{10;\; 50\}$ \\ \hline
\end{tabular}
\caption{Параметры при исследовании зависимости от $\delta$.}
\label{table:2}
\end{table}

В таблице \ref{table:3} указаны параметры задач, на примере которых исследовалась зависимость ошибки от $R$.

\begin{table}[h!]
\centering
\begin{tabular}{|l|c|c|c|c|}
\hline
                        & n                           & $\delta_A$        & $\delta_b$ & $R$         \\ \hline
Шум в матрице           & $10^4$                      & $\{0.01;\; 0.1\}$ &            & $[0;\; 20]$ \\ \hline
Шум в матрице и столбце & \multicolumn{1}{l|}{$10^3$} & $0.001$           & $0.01$     & $[0;\; 50]$ \\ \hline
\end{tabular}
\caption{Параметры при исследовании зависимости от $R$.}
\label{table:3}
\end{table}

Код доступен для просмотра в Google Colab (\href{https://colab.research.google.com/drive/1pZWHk93KahkqwvUrQP8J-Xa1qDEVH78V}{ссылка}).

\paragraph{Результаты}

\begin{enumerate}
    \item Шум в векторе $b$.
    
\begin{figure}[!h]
    \centering
    \subfloat[Враждебный шум в векторе.]{{\includegraphics[width=7.45cm]{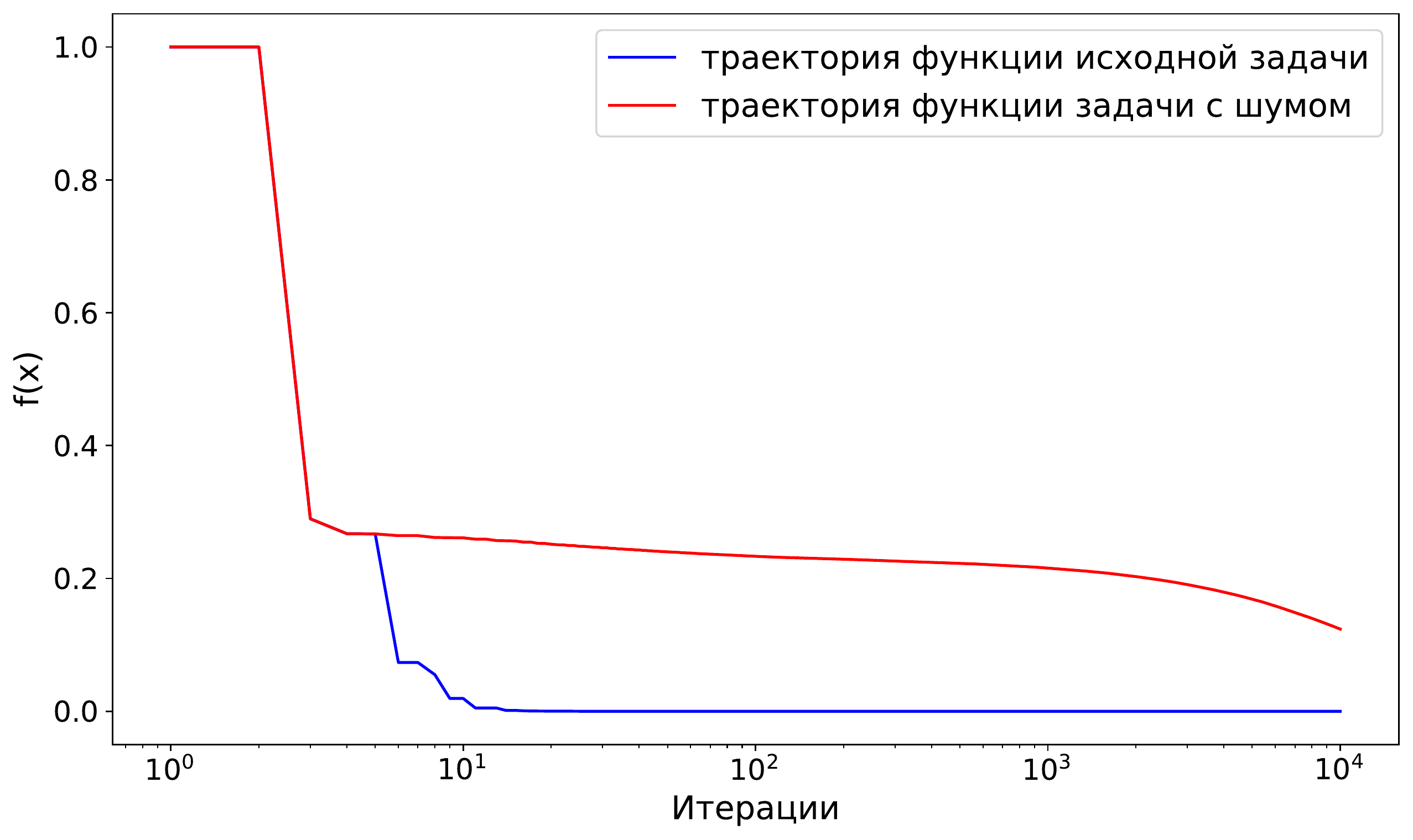} }}%
    \qquad
    \subfloat[Случайный шум в векторе.]{{\includegraphics[width=7.45cm]{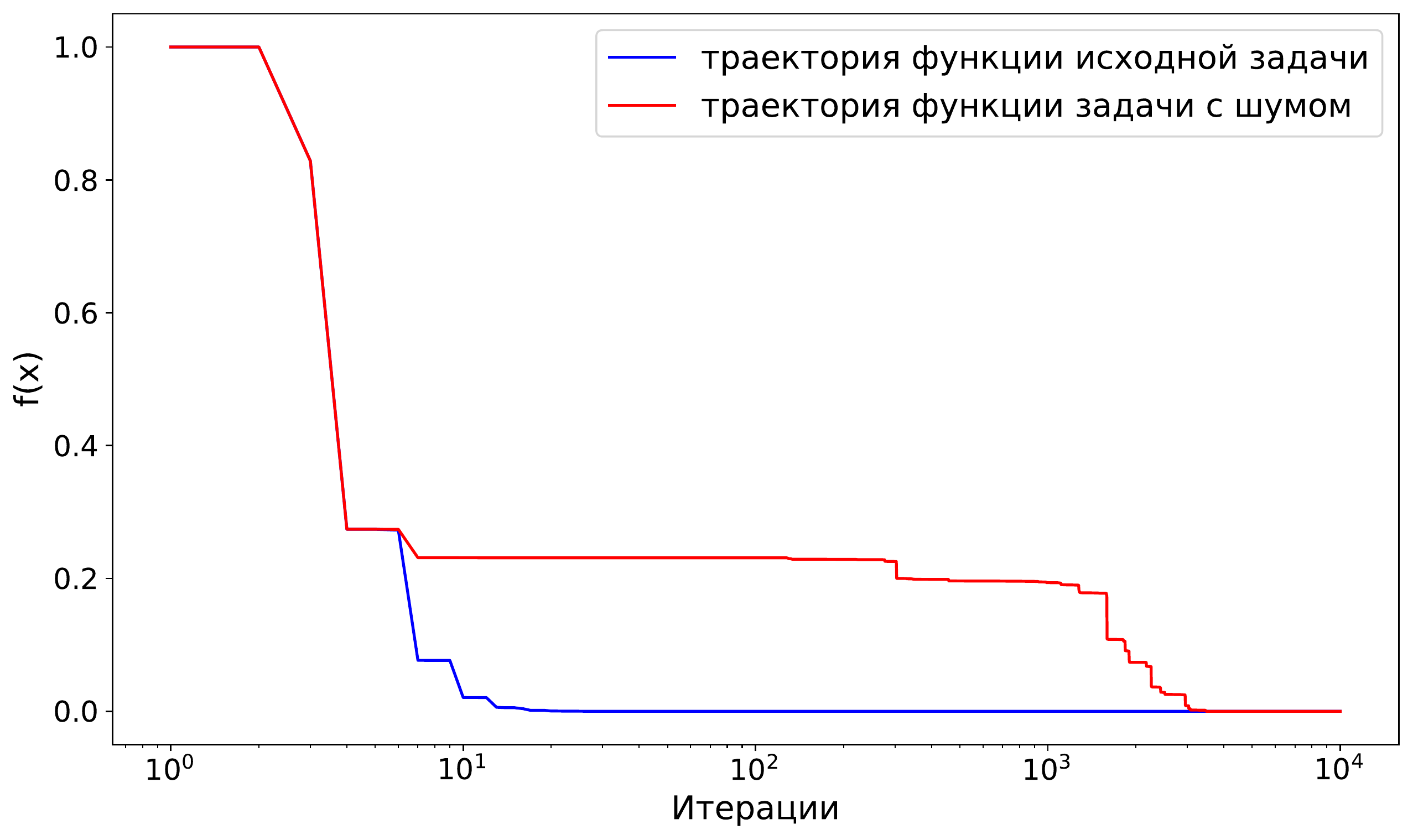} }}%
    \caption{Зависимость величины функции $f(x)$ (отшкалированной на единицу) от номера итерации. Графики показывают выход на асимптоту, свидетельствующий об отсутствии накопления ошибки с ростом числа итераций. Параметры задачи при данном исследовании: $n=10^4,\;\delta_b=0.1,\;R\approx7000$.}%
    \label{fig:2}%
\end{figure}

Накопление ошибки отсутствует как в случае враждебного шума, так и в случае стохастического шума. Метод сходится, но ему требуется больше времени. Этот факт легко установить по графикам на рис. \ref{fig:2}. Из них также видно, что враждебный шум замедляет метод.

\begin{figure}[!h]
    \centering
    \subfloat[Враждебный шум в векторе.]{{\includegraphics[width=7.35cm]{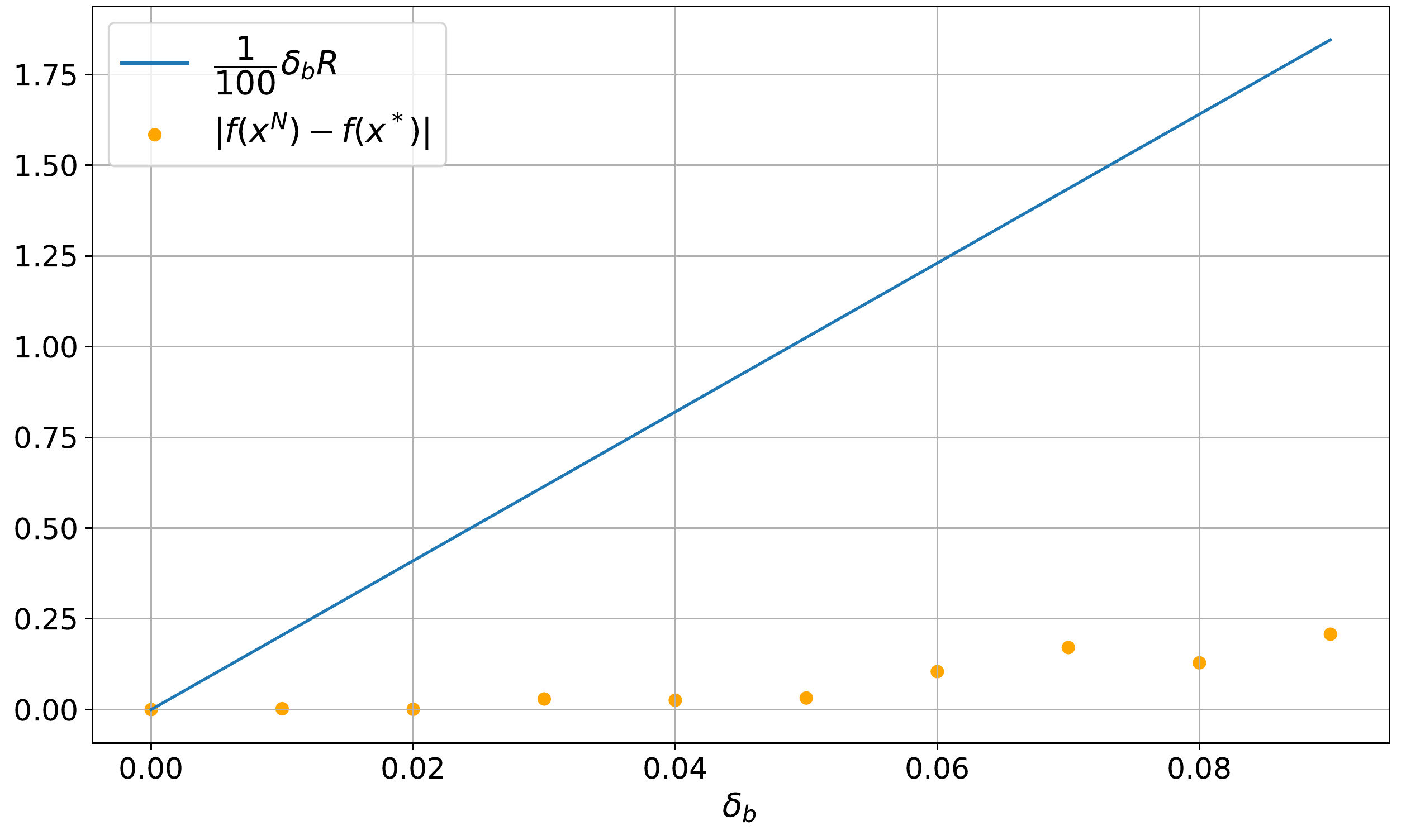} }}%
    \qquad
    \subfloat[Случайный шум в векторе.]{{\includegraphics[width=7.45cm]{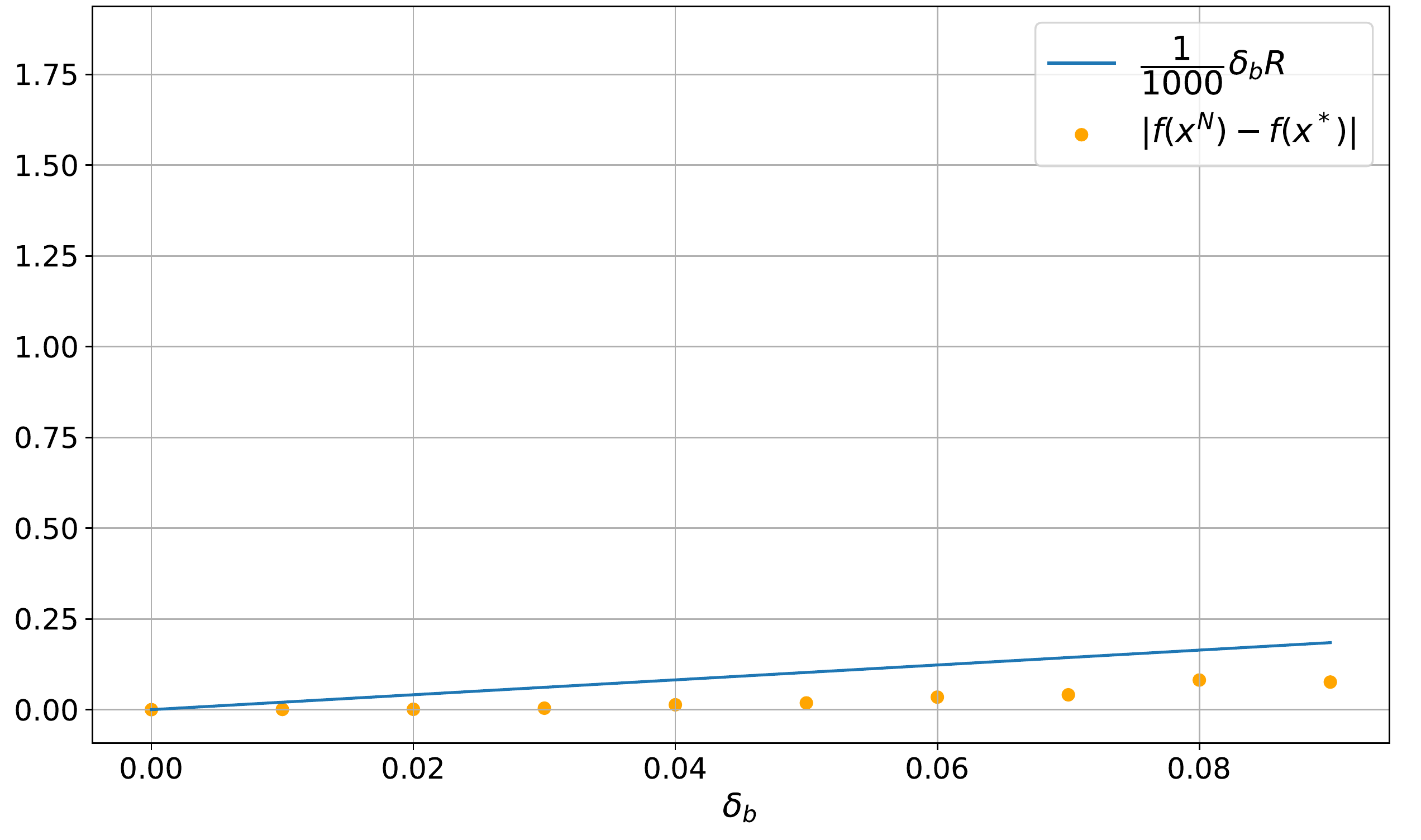} }}%
    \caption{Зависимость невязки по функции от $\delta_b$. Графики показывают, что ошибка зависит от $\delta_b$ линейно. Параметры задачи при данном исследовании: $n=10^4,\; R\approx2000$.}%
    \label{fig:3}%
\end{figure}

Графики на рис. \ref{fig:3} подтверждают гипотезу 
$$|f(x^*_{noisy}) - f(x^*)| = O(\delta_b R)$$ 
при малых\footnote{В данном исследовании имеет смысл говорить исключительно о малых шумах, так как, если шум сравним с нормой градиента, то невозможно найти даже приближённо решение исходной задачи.} $\delta_b$ в случае враждебного и стохастического шума.

\begin{figure}[!h]
    \centering
    \subfloat[Враждебный шум в векторе.]{{\includegraphics[width=7.45cm]{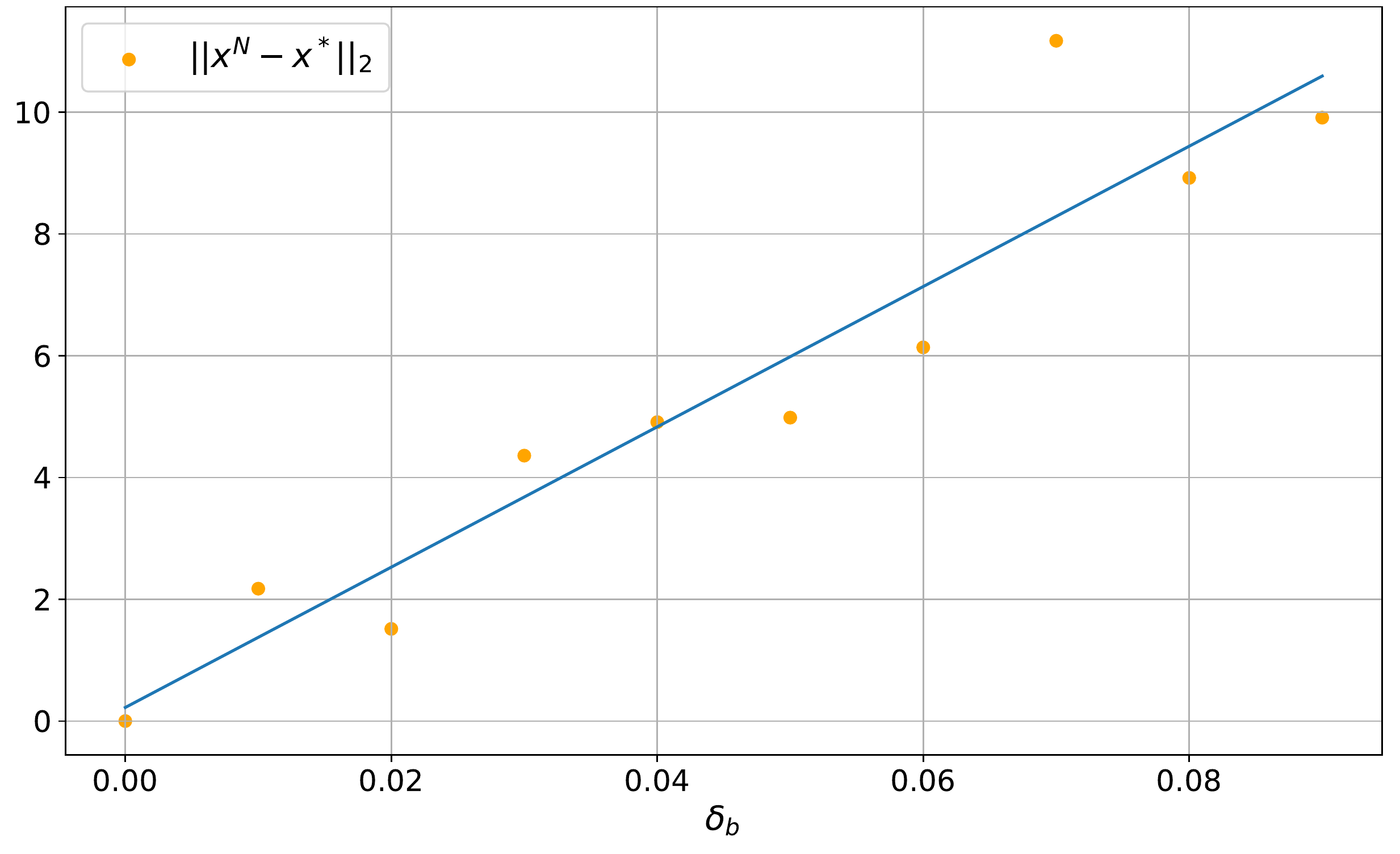} }}%
    \qquad
    \subfloat[Случайный шум в векторе.]{{\includegraphics[width=7.45cm]{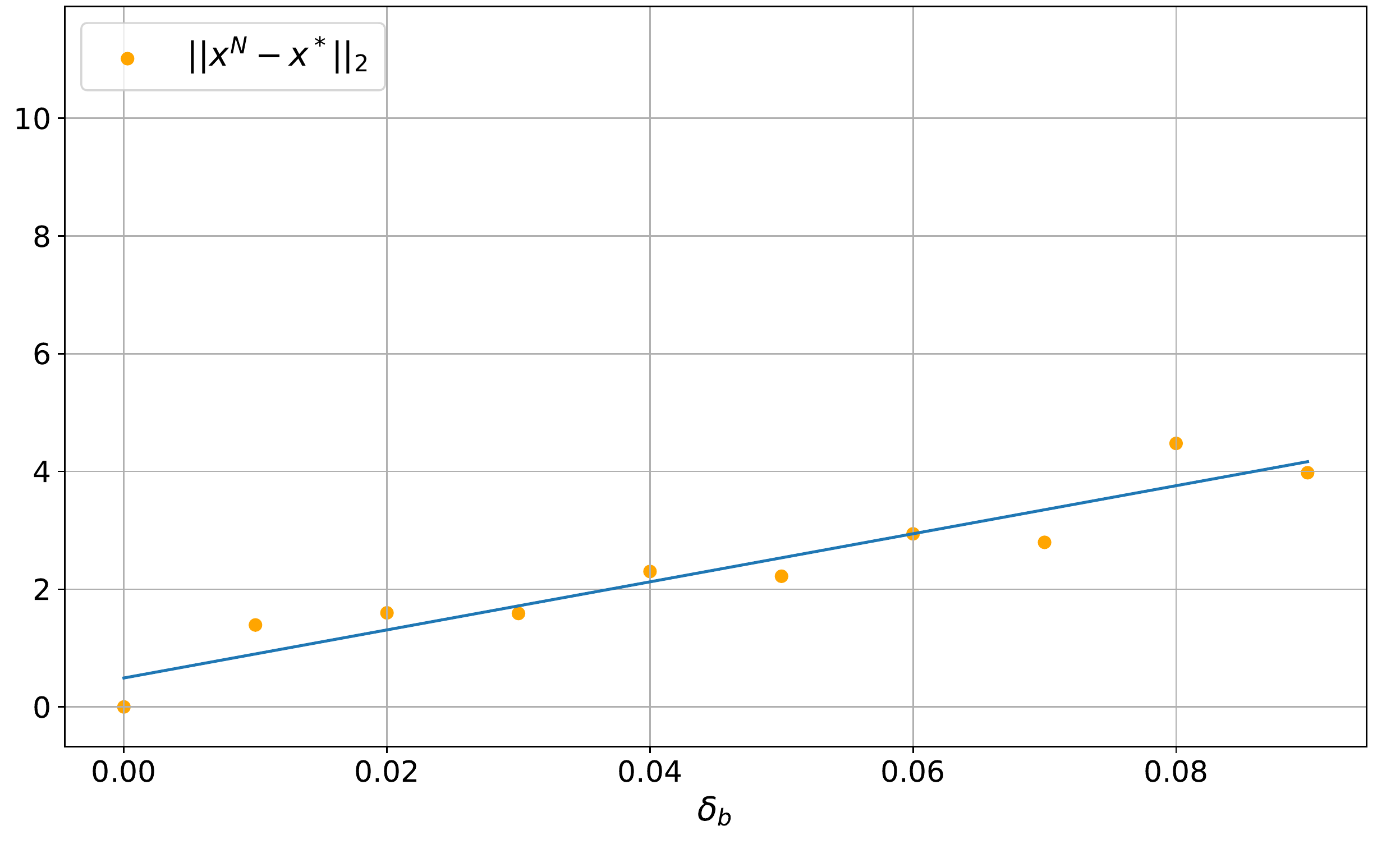} }}%
    \caption{Зависимость невязки по аргументу от $\delta_b$. Параметры задачи: $n=10^4,\; R\approx2000$.}%
    \label{fig:4}%
\end{figure}

Так как задача не является сильно выпуклой, сходимость по аргументу не ожидается. Но тем не менее ошибка аргумента зависит от размера шума линейно, что показано на графиках с рис. \ref{fig:4}. Синяя прямая построена как аппроксимация экспериментальных точек методом наименьших квадратов.

Все перечисленные выше пары графиков чётко показывают, что враждебный шум приводит к более плохим результатам, чем стохастический, этим он обуславливает своё название.

\item Шум в матрице.

\begin{figure}[!h]
    \centering
    \subfloat[Параметры задачи: $n=10^3,\, \delta_A=0.0025,\, R\approx 2000$.]{{\includegraphics[width=7.45cm]{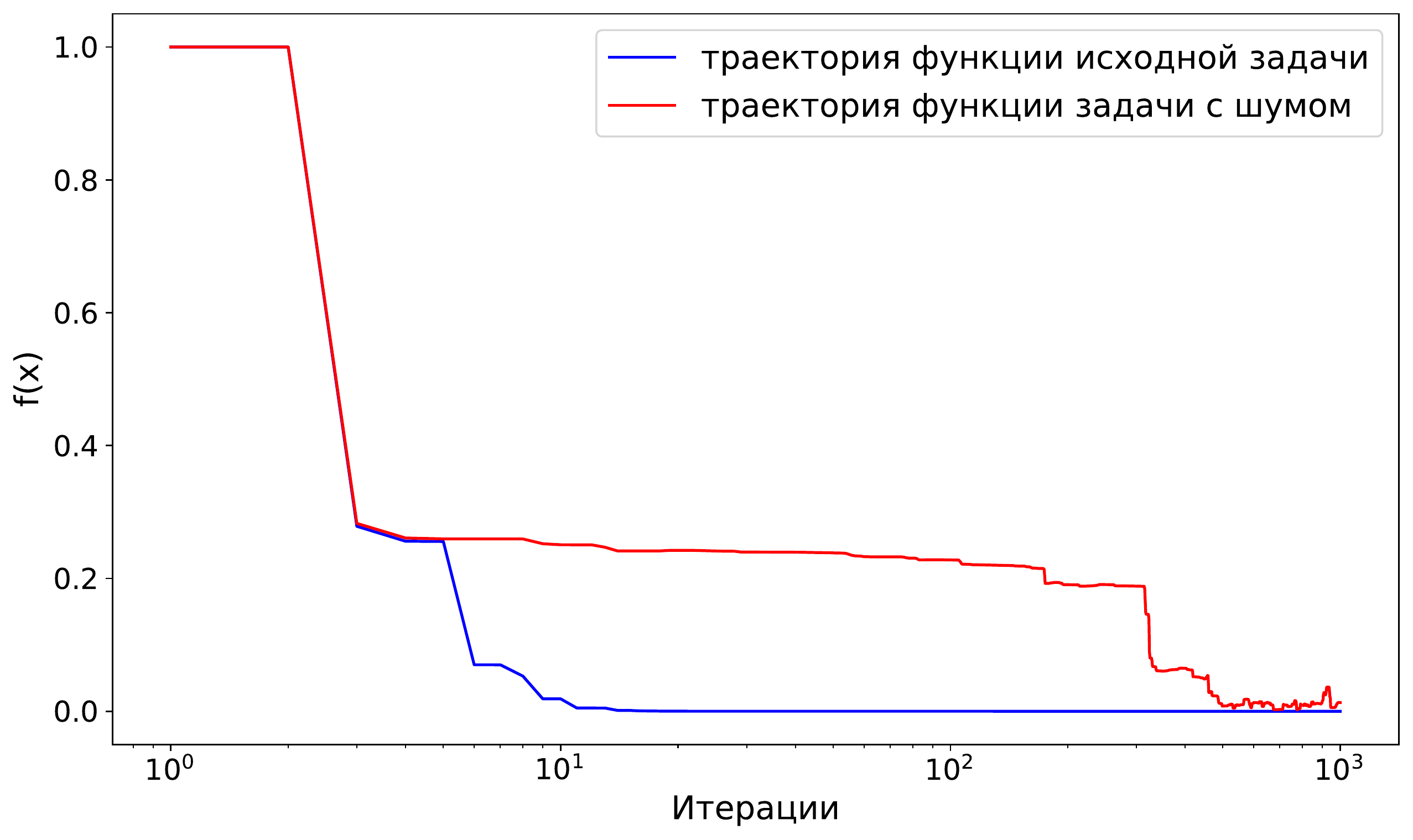} }}%
    \qquad
    \subfloat[Параметры задачи: $n=10^3,\; \delta_A=0.005,\; R\approx 2000$.]{{\includegraphics[width=7.45cm]{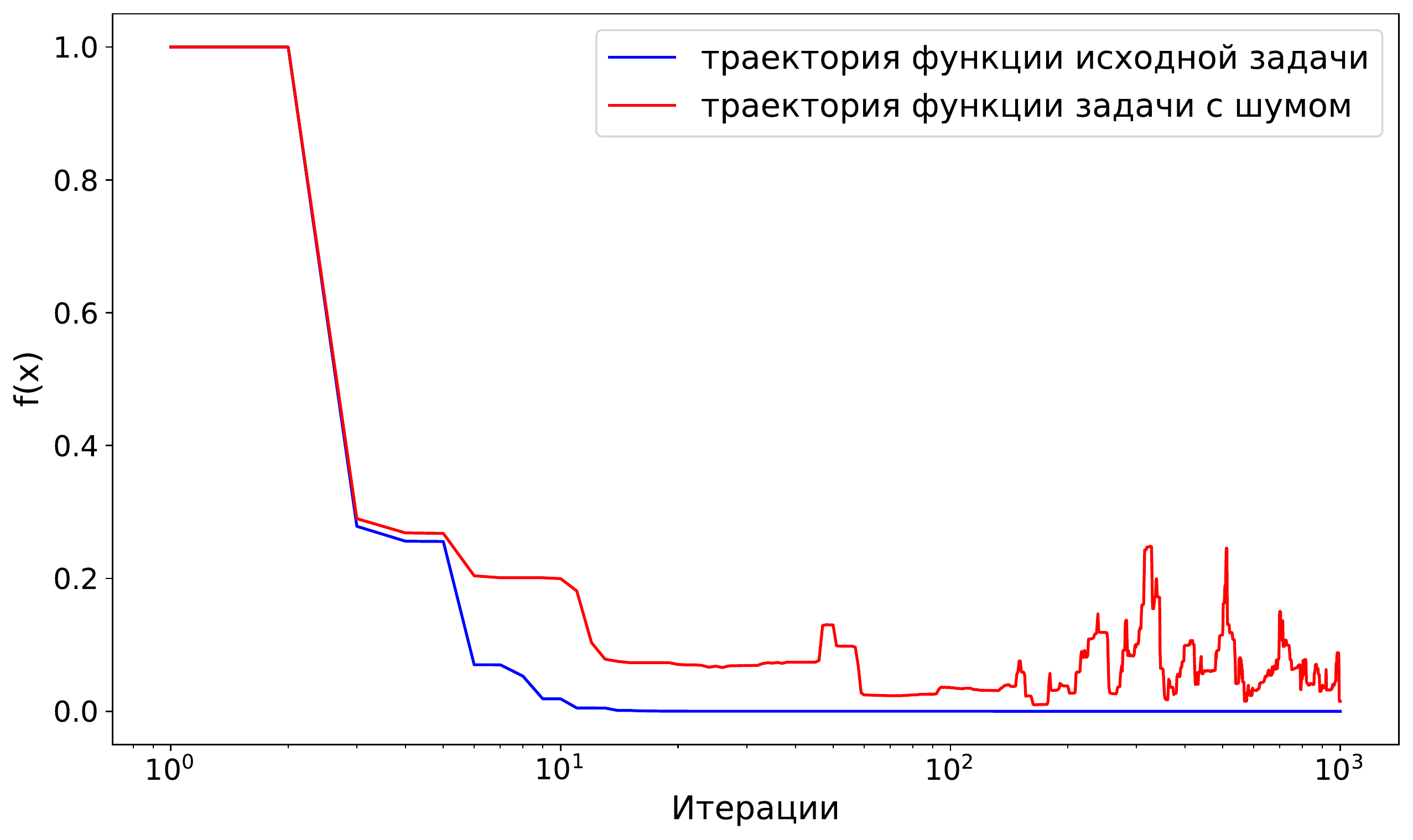} }}%
    \caption{Зависимость величины функции $f(x)$ (отшкалированной на единицу) от номера итерации при шуме в матрице. Графики показывают выход на асимптоту, свидетельствующий об отсутствии накопления ошибки с ростом числа итераций. Также видно, что в случае шума в матрице метод оказывается более чувствительным к величине шума, чем в случае шума в векторе.}%
    \label{fig:5}%
\end{figure}

Если неточность задачи обусловлена шумом в матрице, то накопления ошибки так же не будет (графики на рис. \ref{fig:5}). Но в данном случае метод гораздо более чувствителен к размеру $\delta$. Это обусловлено тем, что шум в градиенте будет ограничен не просто $\delta$, как в случае с шумом в векторе $b$, а $\delta\|x_k\|$, что в абсолютном большинстве итераций (при условии, что $R\gg1$) будет больше, чем просто $\delta$. То есть имея какое-то конкретное $\delta$, ограничение на неточность градиента будет на самом деле больше этого $\delta$.

\begin{figure}[!h]
    \centering
    \subfloat[Параметры задачи: $n=10^4,\; R=50$.]{{\includegraphics[width=7.45cm]{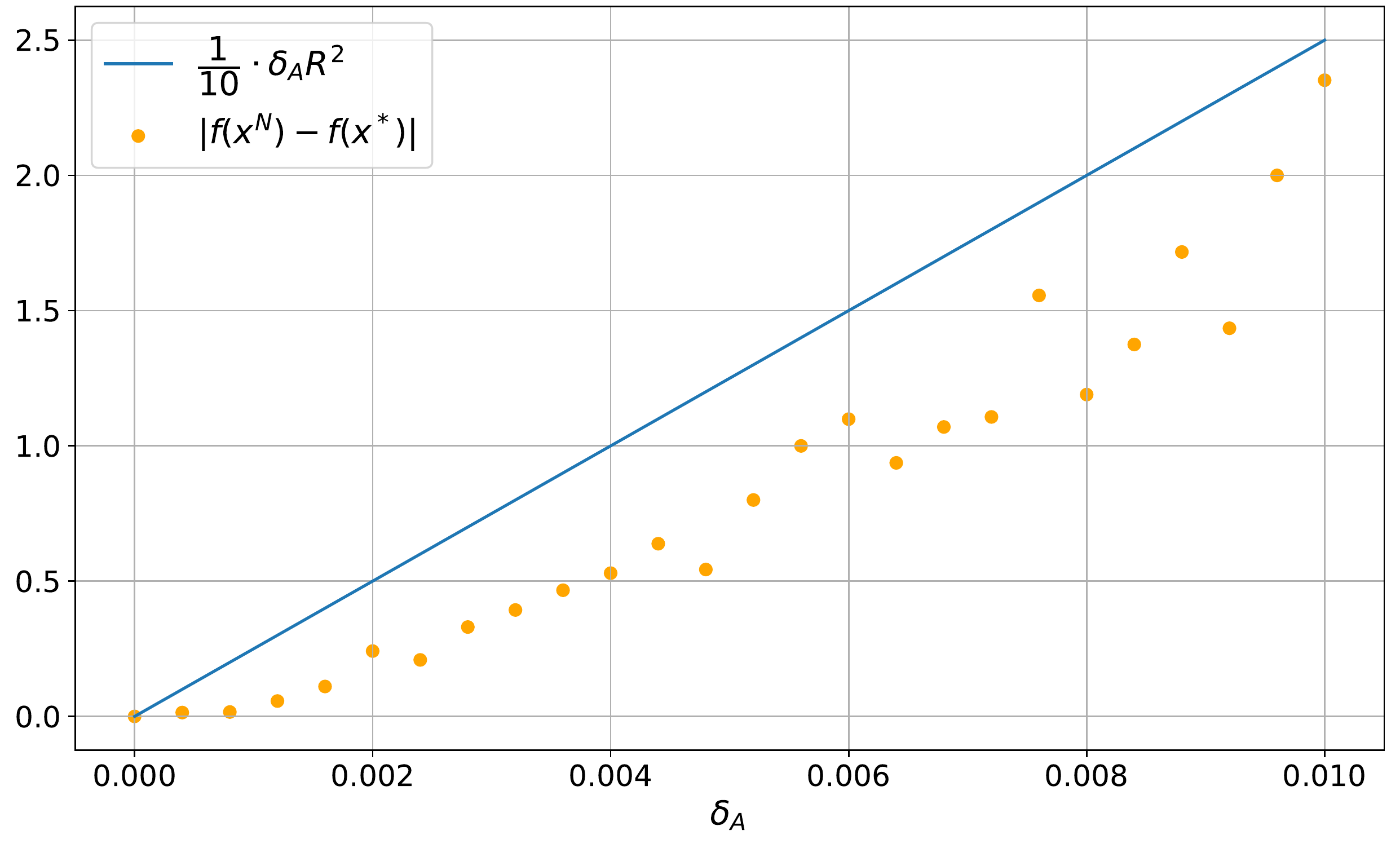} }}%
    \qquad
    \subfloat[Параметры задачи: $n=10^4,\; R=10$.]{{\includegraphics[width=7.45cm]{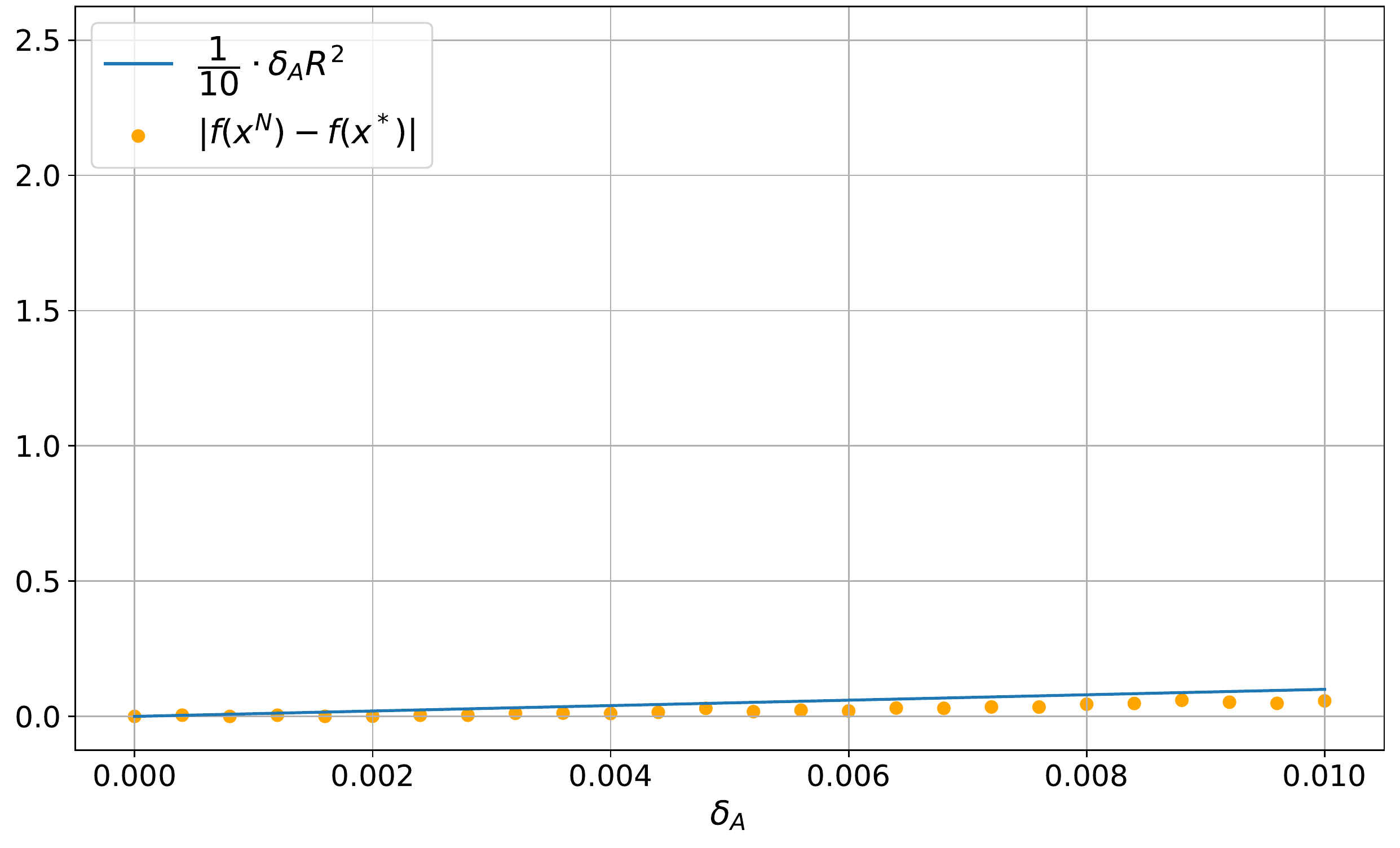} }}%
    \caption{Зависимость невязки по функции от $\delta_A$ при шуме в матрице. Графики показывают, что ошибка зависит от $\delta_A$ линейно. Также видно, что незначительное увеличение $R$ приводит к значительному росту неточности решения по функции.}%
    \label{fig:6}%
\end{figure}

Графики на рис. \ref{fig:6} иллюстрируют линейность зависимости невязки по функции от $\delta_A$, а также сильную зависимость от $R$~---~входящую в оценку, очевидно, не в первой степени. 

\begin{figure}[!h]
    \centering
    \subfloat[Параметры задачи: $n=10^4,\; \delta=0.1$.]{{\includegraphics[width=7.45cm]{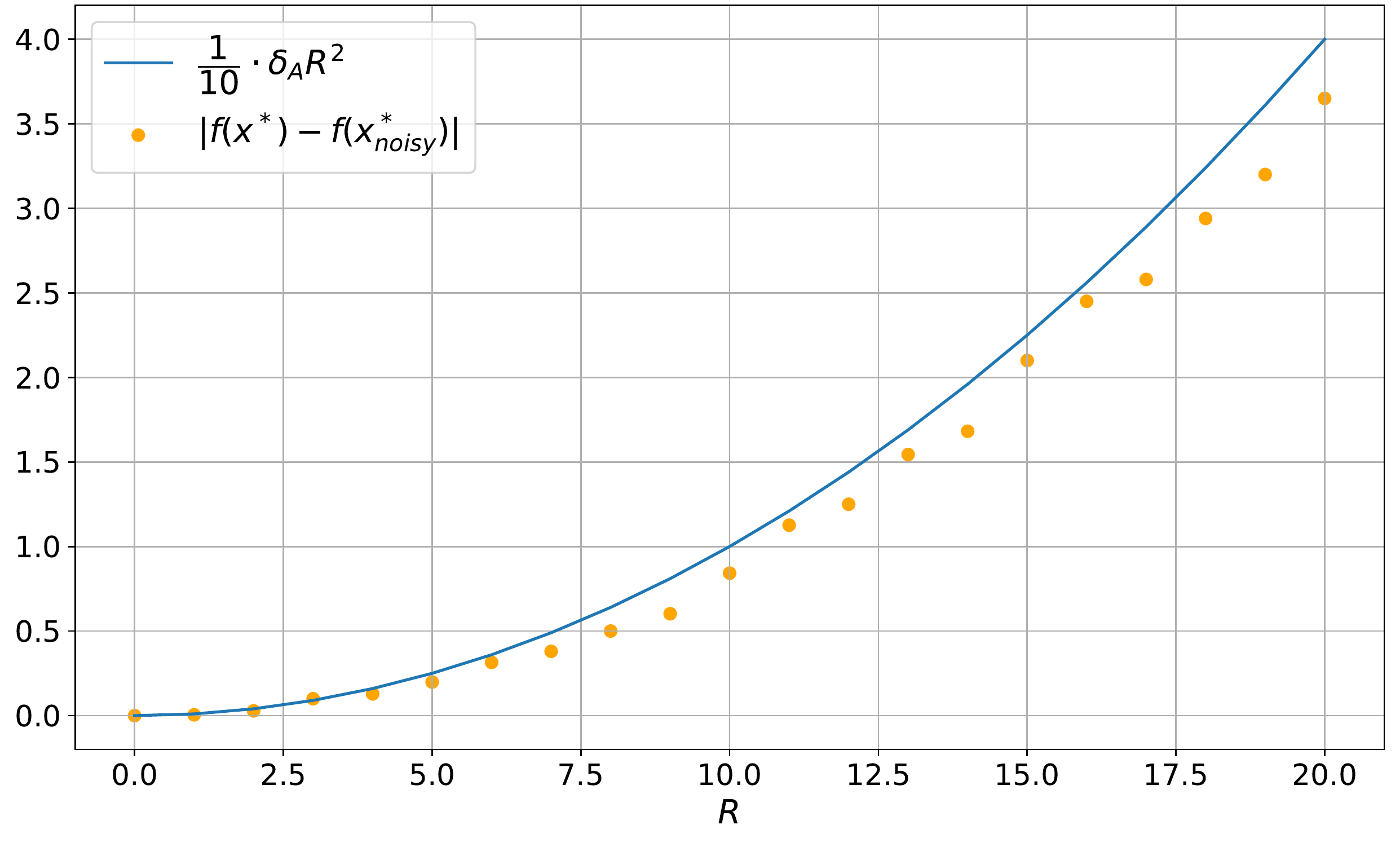} }}%
    \qquad
    \subfloat[Параметры задачи: $n=10^4,\; \delta=0.01$.]{{\includegraphics[width=7.45cm]{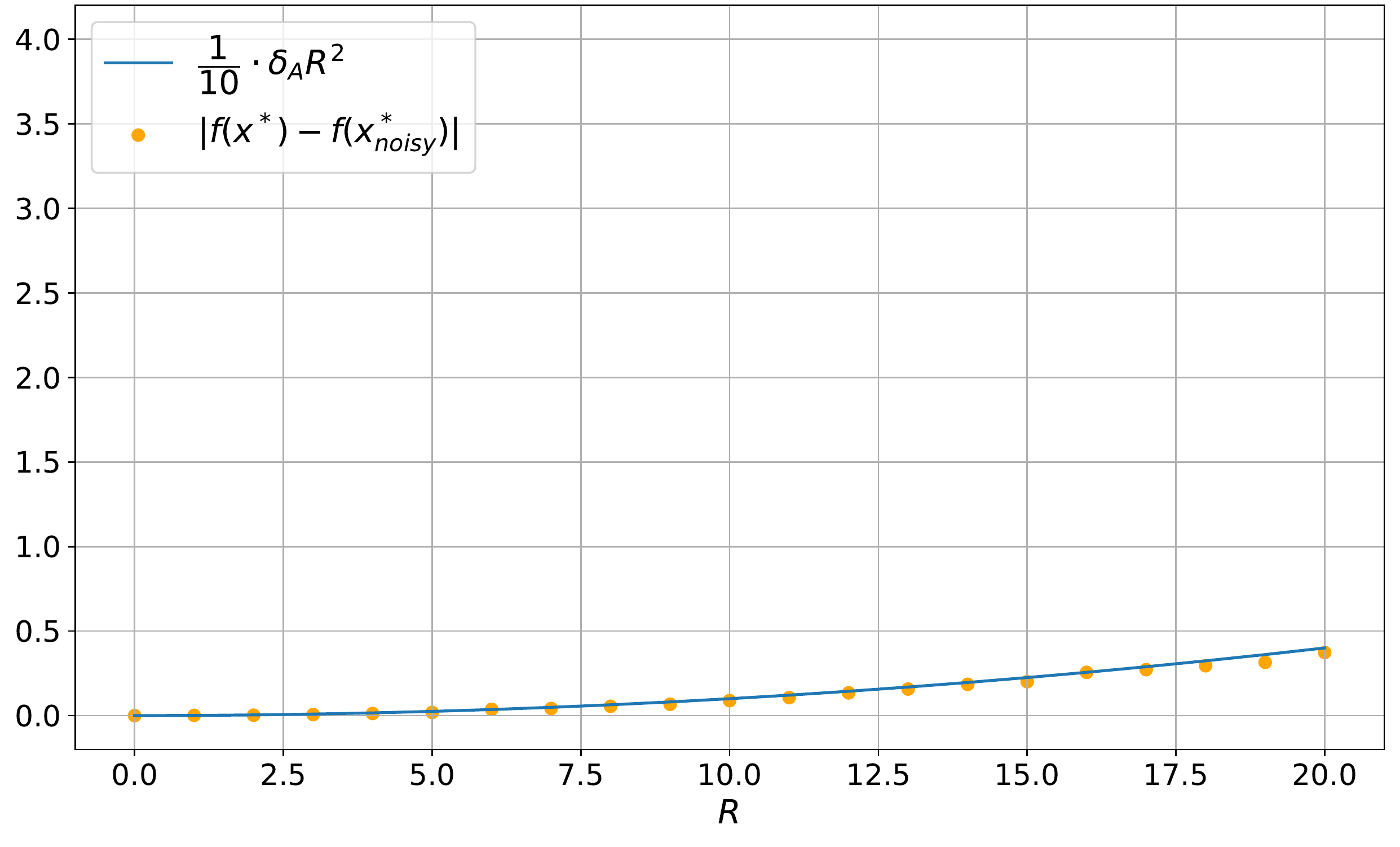} }}%
    \caption{Зависимость невязки по функции от $R$ при шуме в матрице. Графики показывают, что ошибка зависит от $R$ квадратично. Также видно, что увеличение $\delta$ приводит к росту неточности решения по функции.}%
    \label{fig:7}%
\end{figure}

Более ясно иллюстрируют зависимость именно от $R$ графики на рис. \ref{fig:7}. На них изображены результаты исследования зависимости невязки по функции от $R$ при двух различных, но фиксированных размерах $\delta_A$. Эти графики подтверждают гипотезу о вхождении в оценку $$|f(x^*_{noisy}) - f(x^*)| = O(\delta_A R^2)$$ именно квадрата $R$ в отличие от гипотезы $$|f(x^*_{noisy}) - f(x^*)| = O(\delta_b R)$$ для шума в векторе $b$.

Вместе эти четыре графика показывают, что результаты эксперимента согласуются с гипотезой $$|f(x^*_{noisy}) - f(x^*)| = O(\delta_A R^2).$$  

\vspace{0.3cm}

\item Шум в матрице и в векторе.

В данном случае гипотеза является комбинацией гипотез для двух предыдущих вариантов, а именно $$|f(x^*_{noisy}) - f(x^*)| = O(\delta_b R\; +\; \delta_A R^2).$$
Вполне ожидаемо, что и здесь накопление ошибки не наблюдается~---~метод выходит на асимптоту и начинает колебаться около неё, что обусловлено случайностью\footnote{Здесь речь идёт о равновероятном прибавление или вычетание матрицы шума или вектора шума к исходным на каждой итерации. Подробнее в разделе \nameref{supplement materials}.} шума в матрице для графика на рис. \ref{fig:8}(a) и случайностью шумов в матрице и векторе (в большей степени случайностью шума в матрице) для графика на рис. \ref{fig:8}(b). Здесь также видно, что комбинация матричного шума с враждебным шумом в векторе $b$ замедляет метод сильнее, чем комбинация матричного шума со случайным шумом в векторе $b$.

\vspace{0.1cm}

\begin{figure}[!h]
    \centering
    \subfloat[Шум в матрице и враждебный шум в векторе.]{{\includegraphics[width=7.45cm]{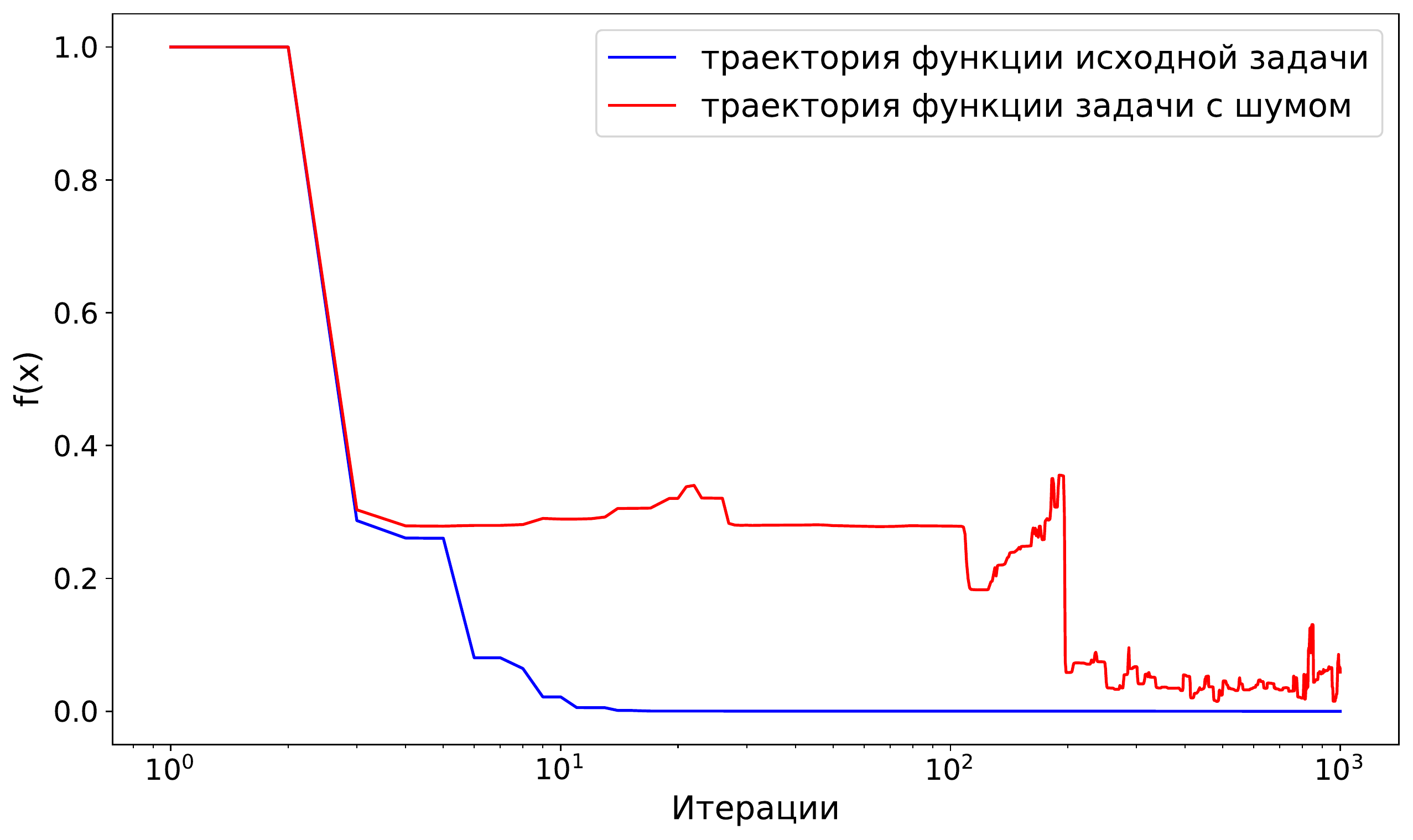} }}%
    \qquad
    \subfloat[Шум в матрице и случайный шум в векторе.]{{\includegraphics[width=7.45cm]{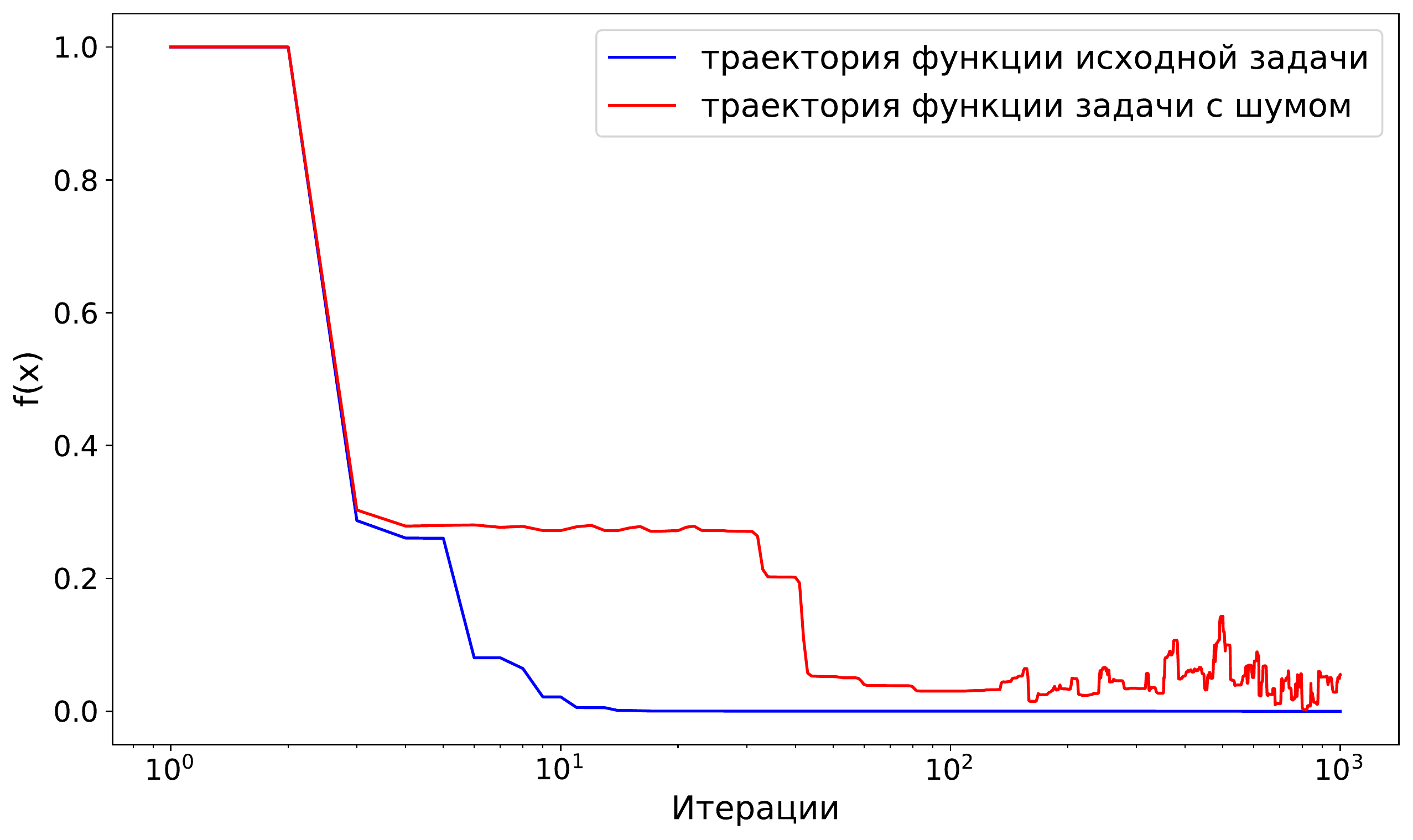} }}%
    \caption{Зависимость величины функции $f(x)$ (отшкалированной на единицу) от номера итерации при шуме в матрице и в векторе. Графики показывают выход на асимптоту, свидетельствующий об отсутствии накопления ошибки с ростом числа итераций. Также видно, что в случае враждебного шума в векторе выход на асимптоту происходит позже, чем в случае случайного шума. Параметры задачи: $n=10^3,\, \delta_A=0.005,\, \delta_b=0.1,\, R\approx2000$.}%
    \label{fig:8}%
\end{figure}

\vspace{0.1cm}

Графики на рис. \ref{fig:9} согласуются с гипотезой. Так же как и предыдущие графики они иллюстрируют, что враждебный шум мешает методу сильнее, чем случайный, и это приводит к большей ошибке по функции найденного решения.

\begin{figure}[!h]
    \centering
    \subfloat[Шум в матрице и враждебный шум в векторе.]{{\includegraphics[width=7.45cm]{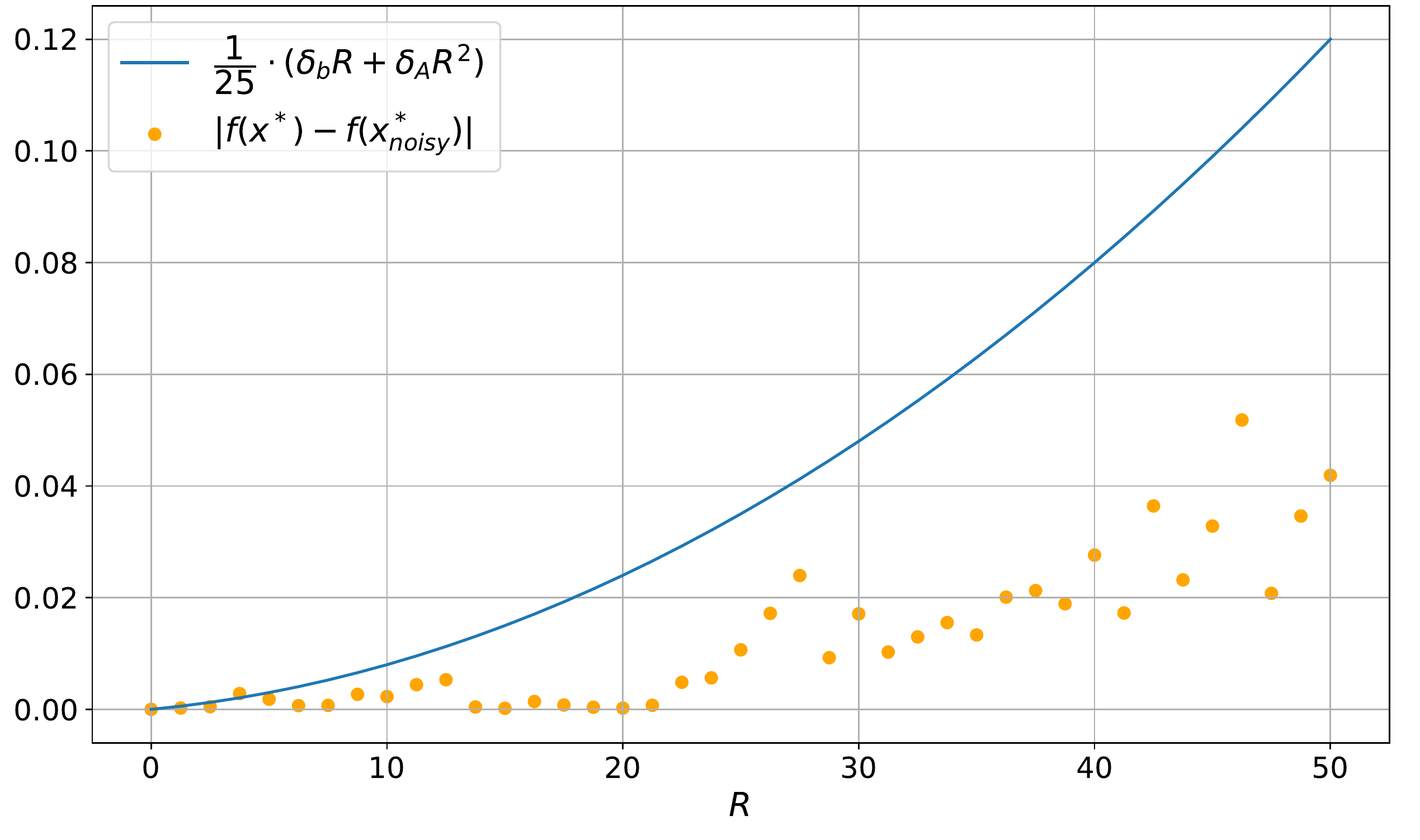} }}%
    \qquad
    \subfloat[Шум в матрице и случайный шум в векторе.]{{\includegraphics[width=7.45cm]{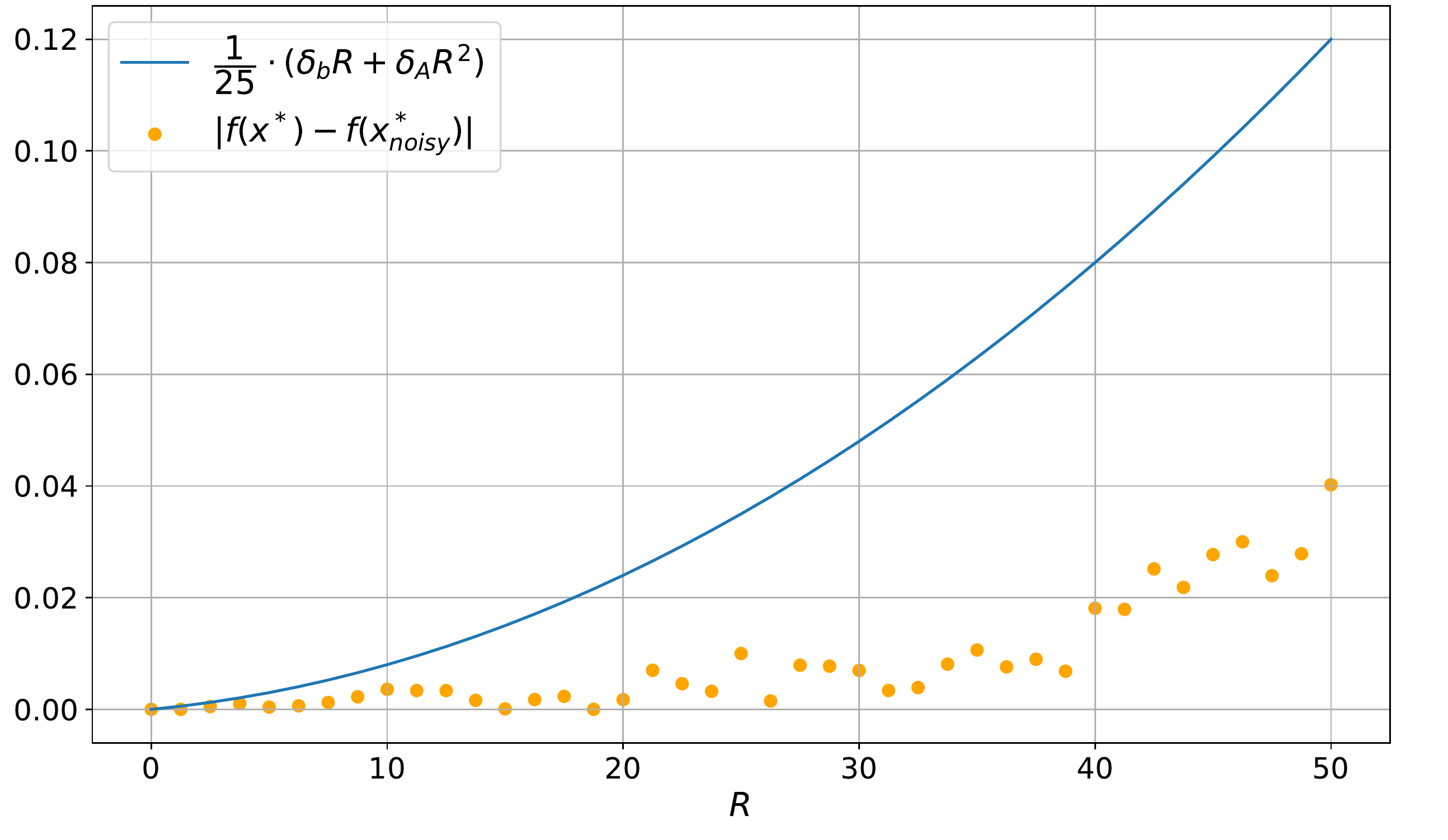} }}%
    \caption{Зависимость невязки по функции от $R$. Графики показывают, что ошибка зависит от $R$ квадратично. Также видно, что шум в матрице вкупе с враждебным шумом в векторе приводит к более плохим результатам, чем шум в матрцие вместе со случайным шумом. Параметры задачи: $n=10^3,\; \delta_A=0.001,\; \delta_b=0.01$.}%
    \label{fig:9}%
\end{figure}
\end{enumerate}

\paragraph{Заключение}

\begin{figure}[!h]
    \centering
    \subfloat[Ускоренный метод Нестерова.]{{\includegraphics[width=7.45cm]{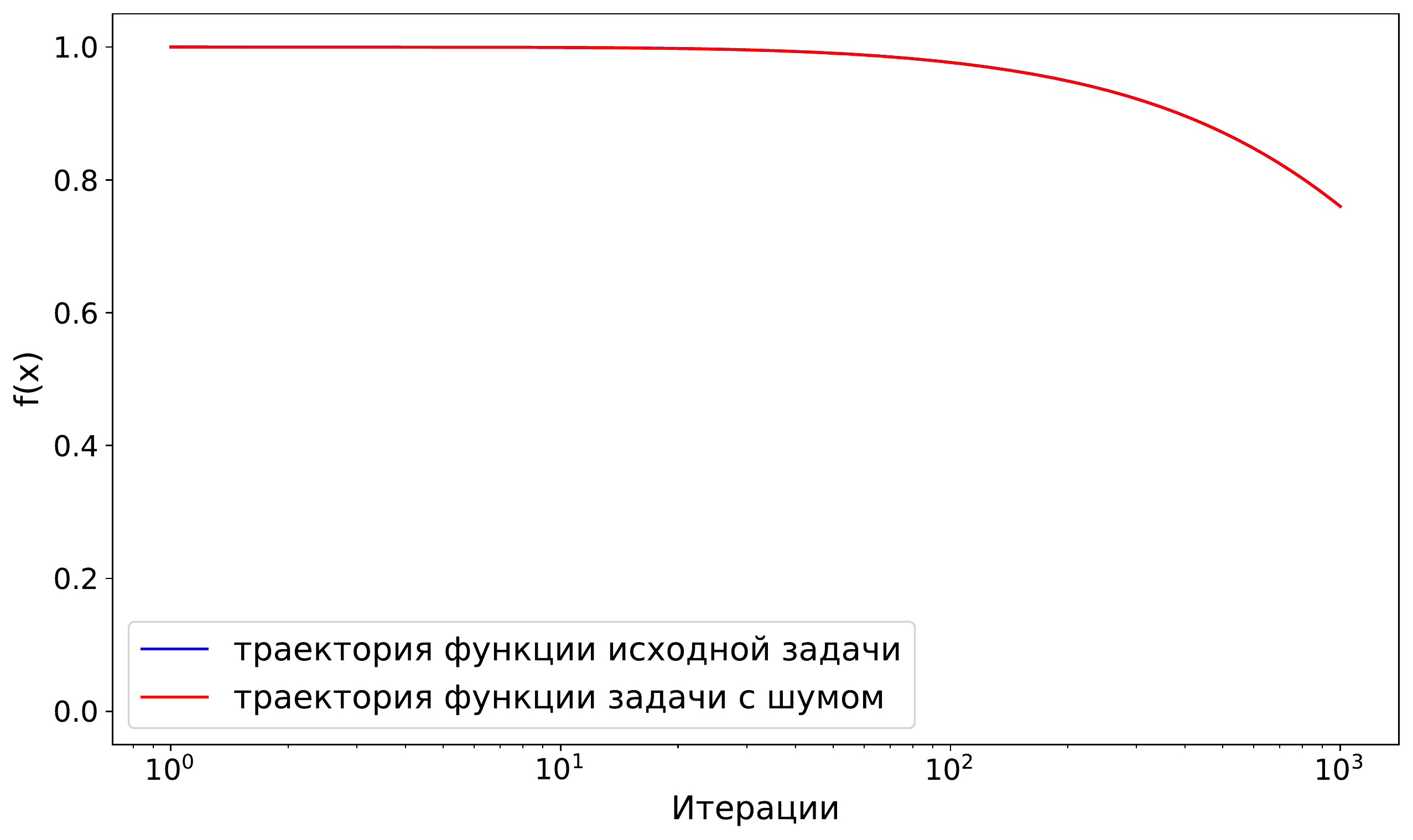} }}%
    \qquad
    \subfloat[Метод сопряжённых градиентов.]{{\includegraphics[width=7.45cm]{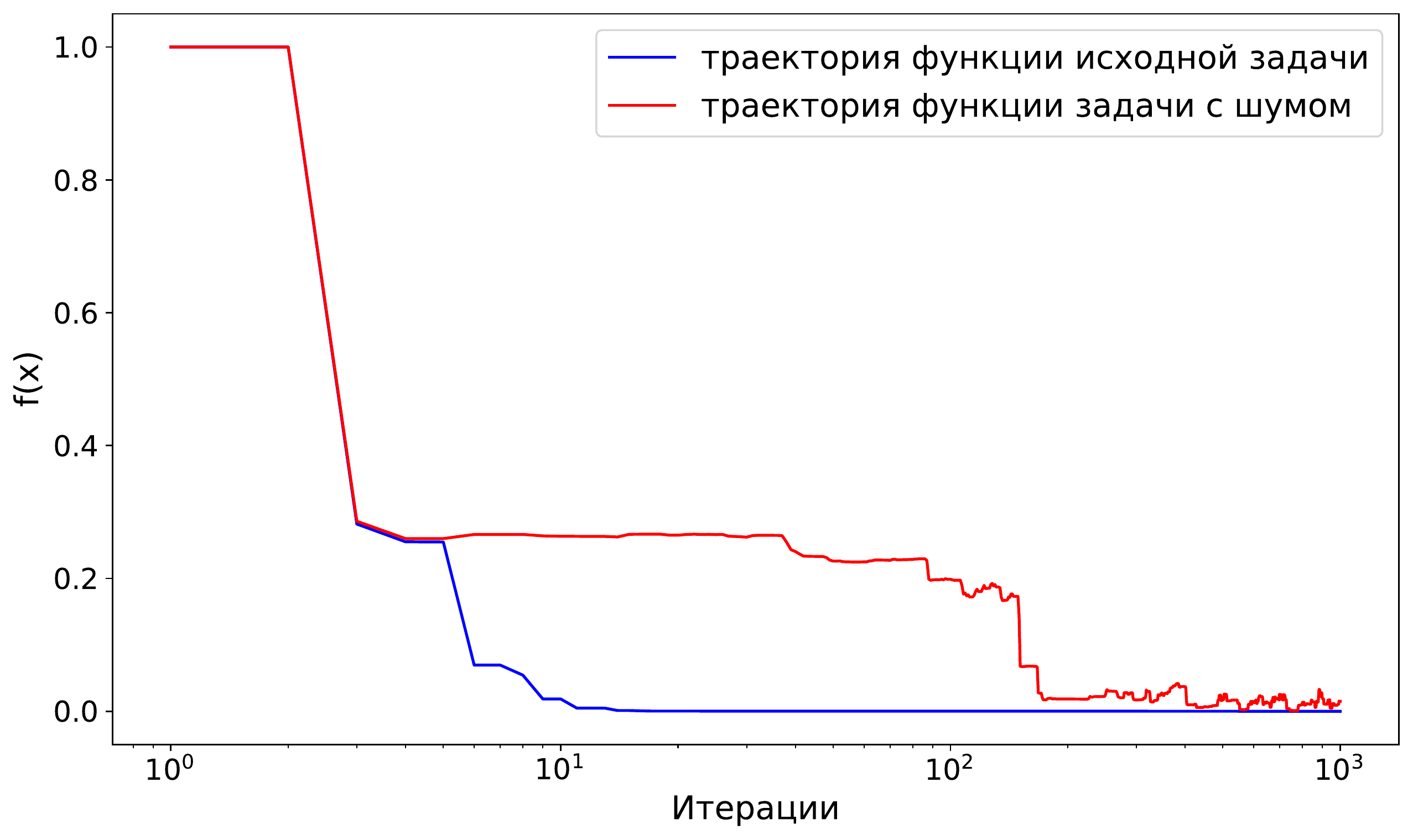} }}%
    \caption{Зависимость величины функции $f(x)$ (отшкалированной на единицу) от номера итерации при шуме в матрице и в векторе. Графики показывают, что метод Нестерова с шумом в оракуле на плохообусловленных задачах сходится сильно медленнее метода сопряжённых градиентов с тем же шумом. Параметры задачи: $n=10^3,\; \delta_A=0.0025,\; \delta_b=0.1$.}%
    \label{fig:10}%
\end{figure}

Результатом данной работы является экспериментальная иллюстрация нетривиального факта: при решении задачи минимизации положительно определённой квадратичной формы методом сопряжённых градиентов с зашумлённым оракулом накопление ошибки отсутствует. При этом, он остаётся наиболее эффективным (быстрым) методом для решения задач минимизации положительно определённых квадратичных форм больших размерностей, что видно из рис. \ref{fig:10}.

Специалистам, сталкивающимся с необходимостью решать СЛАУ или же искать минимумы положительно определённых квадратичных форм с плохообусловленными матрицами больших размерностей, данные экспериментальные результаты добавляют уверенности при использовании метода сопряжённых градиентов в условиях недоступности точного градиента. Хотя шум здесь и приводит к большей неточности полученного решения, чем при использовании других ускоренных методов, метод сопряжённых градиентов оказывается в десятки и даже в сотни раз быстрее. В качестве дальнейшего развития данного исследования было бы интересно обобщить полученный результат на случай гильбертова пространства, где неточность при решении возникает естественным образом ввиду невозможности вычисления градиента по всем направлениям. Также интересно было бы провести исследования скорости сходимости (в вырожденном случае) \textquotedblleft\,в среднем\,\textquotedblright\, по спектру и по выбору точки старта для метода тяжелого шарика, метода Нестерова и метода сопряжённым градиентам \cite{Scieur_Pedregosa20} и стр. 63 \cite{Gasnikov17}.

\paragraph{Дополнительные материалы}
\label{supplement materials}

Важно отметить, что для генерации зашумлённых векторов и матриц есть различные подходы. Выбранные нами способы обусловлены скоростью вычисления и относительной простотой реализации.

\begin{enumerate}
    \item Вектор
        \begin{itemize}
            \item Основа
            \begin{itemize}
                \item $\|{\bf b}\|_2 = \sqrt{b_1^2 + b_2^2 + ... + b_n^2}$
                \item $\|{\bf\tilde{b}} - {\bf b}\|_2 = \sqrt{(\tilde{b}_1 - b_1)^2 + (\tilde{b}_2 - b_2)^2 + ... + (\tilde{b}_n - b_n)^2} \le \delta_b$
        
                \item $\|{\bf\tilde{b}} - {\bf b}\|_2 = \sqrt{\Delta_1^2 + \Delta_2^2 + ... + \Delta_n^2} \le \delta_b$
        
                \item $\Delta_1^2 + \Delta_2^2 + ... + \Delta_n^2 \le \delta_b^2$
            \end{itemize}
            \item Враждебный шум
            \begin{itemize}
                \item $\{\xi_i\}_{i=1}^n \in \mathcal{N}(0, 1)$
        
                \item $\Delta_j^k = \sqrt{\dfrac{\xi_j^2 \cdot \delta_b^2}{\sum_{i=1}^n\xi_i^2}} \cdot \sign\Big(\big[\nabla f(x^k)\big]_j\Big)$, \newline в нашем случае $\sign\Big(\big[\nabla f(x^k)\big]_j\Big) = \sign\Big(\big[Ax^k - b\big]_j\Big)$~--~знак $j$-ого элемента градиента в точке $x^k$
           
                \item ${\bf\Delta^k} = (\Delta_1^k, \Delta_2^k, \dots, \Delta_n^k)^T$
          
                \item ${\bf\tilde{b}^k} = {\bf b} + {\bf\Delta^k}$
            \end{itemize}
            \item Стохастический шум
            \begin{itemize}
                \item $\{\xi_i\}_{i=1}^n \in \mathcal{N}(0, 1)$
        
                \item $\Delta_j^k = \sqrt{\dfrac{\xi_j^2 \cdot \delta_b^2}{\sum_{i=1}^n\xi_i^2}}$~--~нормировка
                
                \item ${\bf\Delta^k} = (\Delta_1^k, \Delta_2^k, \dots, \Delta_n^k)^T$
          
                \item ${\bf\tilde{b}^k} = {\bf b} \pm {\bf\Delta^k}$~--~с вероятностью $\dfrac{1}{2}$
            \end{itemize}
            
        \end{itemize}

    \item Матрица
    \renewcommand{\labelitemi}{--}
    \begin{itemize}
        \item $\|\tilde{A} - A\|_2 \le \delta_A$

        \item $\tilde{A} = A \pm M$~--~с вероятностью $\dfrac{1}{2}$, $\;\; M$~--~матрица шума  

        \item $\tilde{A} - A = \pm M$

        \item $\|M\|_2 \le \delta_A$

        \item $\{\xi_i\}_{i=1}^{n\times n} \in \mathcal{N}(0, 1)$

        \item $m_{1k} = \sqrt{\dfrac{\xi_k^2 \cdot \delta_A^2}{\sum_{i=1}^{n\times n}\xi_i^2}}$~--~нормировка первой строки матрицы шума
        
        \item $m_{2k} = \sqrt{\dfrac{\xi_{(2-1)n + k}^2 \cdot \delta_A^2}{\sum_{i=1}^{n\times n}\xi_i^2}}$~--~нормировка второй строки матрицы шума
        
        \smallskip
        
        \item  \dots
        
        \smallskip
        
        \item $m_{pk} = \sqrt{\dfrac{\xi_{(p-1)n + k}^2 \cdot \delta_A^2}{\sum_{i=1}^{n\times n}\xi_i^2}}$~--~нормировка $p$-ой строки матрицы шума
        
        \smallskip
        
        \item  \dots
        
        \item $m_{nk} = \sqrt{\dfrac{\xi_{(n-1)n + k}^2 \cdot \delta_A^2}{\sum_{i=1}^{n\times n}\xi_i^2}}$~--~нормировка последней строки матрицы шума

        
        \item $M = \begin{bmatrix}
        m_{11} & m_{12} & \dots & m_{1n} \\ 
        m_{21} & m_{22} & \dots & m_{2n} \\ 
        \vdots & \vdots & \ddots & \vdots \\ 
        m_{n1} & m_{n_2} & \cdots & m_{nn}\end{bmatrix}$
    \end{itemize}

\end{enumerate}



\begin{thebibliography}{99}

\bibitem[Gasnikov, 2017]{Gasnikov17}
	\textit{Gasnikov~A.\,V.} Universal gradient descent //arXiv preprint arXiv:1711.00394. – 2017.

\bibitem[Poljak, 1981]{poljak81}
	\textit{Poljak~B.\,T.} Iterative algorithms for singular minimization problems //Nonlinear Programming 4. – Academic Press, 1981. – С. 147-166.
	
\bibitem[Нестеров, 2010]{nesterov10}
	\textit{Нестеров~Ю.\,Е.} Введение в выпуклую оптимизацию. – 2010.
	
	\vspace{0.1cm}{\footnotesize{\it Nesterov~Yu.\,E.} Vvedenie v vypukluyu optimizatsiyu. – 2010.\par}
	
\bibitem[Гилл, Мюррей, РАйт, 1985]{gillmyurreirait85}
	\textit{Гилл~Ф., Мюррей~У., Райт~М.} Практическая оптимизация. – Мир, 1985. – Т. 509.
	
	\vspace{0.1cm}{\footnotesize{\it Gill~F., Myurrei~U., Rait~M.} Prakticheskaya optimizatsiya. – Mir, 1985. – T. 509.\par}

\bibitem[Немировский, 1986]{nemirovskiy86}
	\textit{Немировский~А.\,C.} О регуляризующих свойствах метода сопряжённых градиентов на некорректных задачах //Журнал вычислительной математики и математической физики. – 1986. – Т. 26. – №. 3. – С. 332-347.
	
	\vspace{0.1cm}{\footnotesize{\it Nemirovskii~A.\,S.} O regulyarizuyushchikh svoistvakh metoda sopryazhennykh gradientov na nekorrektnykh zadachakh //Zhurnal vychislitel'noi matematiki i matematicheskoi fiziki. – 1986. – T. 26. – №. 3. – S. 332-347.\par}

\bibitem[Немировский А.С., Поляк Б.Т., 1984]{nemirovskiy84}
	\textit{Немировский~А.\,C., Поляк~Б.\,Т.} Итерационные методы решения линейных некорректных задач при точной информации. II. // Изв. АН СССР. Техническая кибернетика – 1984 – № 3. – С. 18–25.
	
	\vspace{0.1cm}{\footnotesize{\it Nemirovskii~A.\,S., Polyak~B.\,T.} Iteratsionnye metody resheniya lineinykh nekorrektnykh zadach pri tochnoi informatsii. II. // Izv. AN SSSR. Tekhnicheskaya kibernetika – 1984 – № 3. – S. 18–25.\par}

\bibitem[Nemirovski A., 1992]{nemirovskiy92}
	\textit{Nemirovski~A.} Information-based complexity of linear operator equations // Journal of Complexity. – 1992. – V. 8. – P. 153–175.
	
\bibitem[Kabanikhin S.I., 2012]{kabanikhin2012}
	\textit{Kabanikhin~S.\,I.} Inverse and ill-posed problems. – De Gruyter, 2012.

\bibitem[Devolder O., 2013]{devolder2013}
	\textit{Devolder~O.} Exactness, inexactness and stochasticity in first-order methods for large-scale convex optimization: PhD thesis. – CORE UCL, March 2013.

\bibitem[Dvinskikh D., Gasnikov A., 2019]{dvinskikh2019}
	\textit{Dvinskikh~D., Gasnikov~A.} Decentralized and parallelized primal and dual accelerated methods for stochastic convex programming problems //arXiv preprint arXiv:1904.09015. – 2019. 
 
\bibitem[d'Aspremont A., 2008]{dAspremont}
	\textit{d'Aspremont~A.} Smooth optimization with approximate gradient //SIAM Journal on Optimization. – 2008. – Т. 19. – №. 3. – С. 1171-1183.
	
	
\bibitem[Scieur D., Pedregosa F., 2020]{Scieur_Pedregosa20}
	\textit{Scieur~D., Pedregosa~F.} Universal Average-Case Optimality of Polyak Momentum //arXiv preprint arXiv:2002.04664. – 2020.
 



\end{thebibliography}
\end{document}